\documentclass[10pt]{article}
\usepackage{}
\usepackage{amsfonts}
\usepackage{mathrsfs}
\input{epsf.sty}
\usepackage{epsfig}
\usepackage{amsmath}
\usepackage{amssymb}
\usepackage{amsthm}
\usepackage{times}
\usepackage{float}
\usepackage{color}

\usepackage{indentfirst}
\def\be{\begin{equation}}
\def\ee{\end{equation}}
\def\ba{\begin{array}}
\def\ea{\end{array}}

\newtheorem{thm}{Theorem}[section]

\newtheorem{lem}[thm]{Lemma}
\newtheorem{prop}[thm]{Proposition}

\newtheorem{rem}{Remark}

\newcommand{\mi}{\mbox{i}}
\numberwithin{equation}{section}%\numberwithin{equation}{section}

\newcommand{\norm}[1]{\left\Vert#1\right\Vert}
\newcommand{\abs}[1]{\left\vert#1\right\vert}

\newcommand{\md}{\mathrm{d}}

\def\be{\begin{equation}}
\def\ee{\end{equation}}
\def\br{\begin{eqnarray}}
\def\er{\end{eqnarray}}

\begin{document}
\title{ Quasi-periodic solutions for perturbed generalized nonlinear vibrating string equation with singularities \footnote{Supported by NNSFC 11271076 and NNSFC 11121101.} }

\author{ $\mbox{Chengming \ Cao}^{\ddag}$ \hspace{12pt}  \  \  Xiaoping Yuan$^{\dag}$\\
${}^{\ddag}\mbox{School of Mathematical Sciences, Fudan University,
Shanghai 200433, P R China}$\\$ \mbox {Email:12110180001@fudan.edu.cn}$\\
$^\dag$ School of Mathematical Sciences, Fudan University,
Shanghai 200433, P R China\\
Email: xpyuan@fudan.edu.cn}

\date{}
\maketitle
 {\bf Abstract}: The existence of 2-dimensional KAM tori is proved for the perturbed  generalized nonlinear vibrating string equation with singularities $u_{tt}=((1-x^2)u_x)_x-mu-u^3$ subject to certain boundary conditions by means of infinite-dimensional KAM theory with the help of partial Birkhoff normal form, the characterization of the singular function space and the estimate of the integrals related to Legendre basis.

{\bf Key words:}  KAM theory; quasi-periodic solutions; singular differential operator

\section{Introduction}

The KAM (Kolmogorov-Arnold-Moser) theory was used to find the quasi-periodic solutions for hamiltonian partial differential equations (PDEs), originally  by Kuksin \cite{§Ü§å§Ü§ã§Ú§ßnearly,kuksin1987hamiltonian,kuksin1989perturbation}  and Wayne  \cite{wayne1990periodic}. Among those PDEs, the nonlinear wave (NLW) equation \begin{equation}\label{nlw}
  u_{tt}-u_{xx}+Vu+f(u)=0, f(u)=\sum_{k\geq 3}f_ku^k
\end{equation}
has been investigated by many authors.

In KAM theory some parameters are needed to overcome resonances arising in the small divisors.  Kuksin \cite{§Ü§å§Ü§ã§Ú§ßnearly} assumed that the  potential $V=V(x;\xi)$ depends on an $n$-dimensional parameter vector $\xi$ and   showed that there are many quasi-periodic solutions for NLW for ``most'' parameters $\xi'$s. See also \cite{bourgain1994construction,wayne1990periodic,grebert2011kam}. See  P\"{o}schel \cite{poschel1996quasi} for  constant-value potential $V(x)\equiv m$ with $m>0$ and $-1<m<0$  and   \cite{yuan2003invariant} for $V(x)\equiv m\in(-\infty,-1)\setminus \mathbb{Z}$ and  \cite{yuan2006invariant} for any prescribed nonconstant potential $V\in L^2[0,\pi]$.  When $V(x)\equiv 0$ which is called completely resonant, Berti and Procesi \cite{berti2006quasi} proved the existence of 2-dimensional tori and the existence of any dimensional KAM tori was proved in \cite{yuan2006quasi} .

In the above papers, the potentials $V$ are regular. In physics and mechanics the potentials sometimes contain some kind of singularity.  As an example, let us consider the Legendre potential,
\begin{equation}
  V_L(x)=-\frac{1}{2}-\frac{1}{4}\tan^2{x},\quad x\in[-\frac{\pi}{2},\frac{\pi}{2}].
\end{equation}
Since
\begin{equation*}
\lim_{x\rightarrow\pm \frac{\pi}{2}}V_L(x)=-\infty,
\end{equation*}
the endpoints $x=\pm \frac{\pi}{2}$ are actually singular.

It is well-known that the singular differential expression
\begin{equation}\label{A1}
\tilde{\mathscr{A}}:=-\frac{\md^2}{\md x^2}+V_L(x),\quad x\in [-\frac{\pi}{2},\frac{\pi}{2}],
\end{equation}
is in limit-circle case and is of deficiency index $(2,2)$. The expression $\tilde{\mathcal{A}}$ is a self-adjoint operator in the domain
\begin{equation*}
\mathscr{D}(\tilde{\mathscr{A}})=\left\{u(x)\in L^2\left[-\frac{\pi}{2},\frac{\pi}{2}\right]\quad\Big{|}\: u\left(\pm \frac{\pi}{2}\right)=0\right\}.
\end{equation*}
Introducing the change of variable
\begin{equation*}\left\{
\begin{aligned}
&y=\sin{x},\quad x\in\left[-\frac{\pi}{2},\frac{\pi}{2}\right],\\
&z=\frac{u}{\sqrt{\cos{x}}},
\end{aligned}\right.
\end{equation*}
the operator $\tilde{\mathscr{A}}$ with its domain can be written as
\begin{equation}
  \mathscr{A}= -\frac{\md}{\md y}(1-y^2)\frac{\md}{\md y},\quad y\in[-1,1]
\end{equation}
with
\begin{equation}\label{space}
\begin{aligned}
\mathscr{D}(\mathscr{A})=&\left\{z(y)\in L^2[-1,1]\Big{|} \lim_{y\rightarrow \pm 1}z(y)\sqrt{\cos{\arcsin{y}}}=0\right\}\\
=&\left\{z(y)\in L^2[-1,1]\Big{|} \lim_{y\rightarrow \pm 1}z(y)(1-y^2)^{\frac{1}{4}}=0\right\},
\end{aligned}
\end{equation}
In convention, we still write $z(y)=u(x), y=x$.
The operator $\mathcal{A}$ has pure point spectrum $\sigma(\mathscr{A})=\sigma_p(\mathscr{A})$. And the property
\begin{equation*}
(\mathscr{A}u,u)=\int_{-1}^{1}-\frac{\md }{\md x}\Big[(1-x^2)\frac{\md}{\md x} u(x)\Big]\overline{u(x)}\md x =\int_{-1}^{1}(1-x^2)\Big|\frac{\md}{\md x}u(x)\Big|^2\md x\geq 0
\end{equation*} yields
\begin{equation*}
\sigma(\mathscr{A})\subset[0,\infty).
\end{equation*}

To ensure the singular differential operator's strict positive definiteness, we use the notation
\begin{equation}
  A=\mathscr{A}+m \quad (m>0).
\end{equation}
Let $\lambda^2_j$ and $\phi_j(j=1,2,\ldots)$ be the eigenvalues and eigenfunctions of $A$,respectively. Here $\lambda_j>0 (j=1,2,\ldots)$.

Write
\begin{equation}\label{u}
u(t,x)=\sum_{j\geq1}\frac{q_j(t)}{\sqrt{\lambda_j}}\phi_j(x)
\end{equation}
Inserting (\ref{u}) into the following equation
\begin{equation*}
\left\{
\begin{aligned}
&u_{tt}-\big((1-x^2)u_x\big)_x+mu+ u^3=0,\quad x\in[-1,1].\\
&\lim_{x\rightarrow\pm 1}u(x)(1-x^2)^{\frac{1}{4}}=0
\end{aligned}
\right.
\end{equation*}
we have
\begin{equation}
\ddot{q_j}+\lambda_j^2q_j+\sqrt{\lambda_j} \langle u^3,\phi_j\rangle=0.
\end{equation}
This is a hamiltonian system
\begin{equation}\label{qp}
\left\{
  \begin{aligned}
  &\dot{q_j}=\frac{\partial H}{\partial p_j}=\lambda_jp_j\\
  &\dot{p_j}=-\frac{\partial H}{\partial q_j}=-\lambda_jq_j-\frac{\partial G}{\partial q_j}
  \end{aligned}\quad,
\right.\quad j=1,2,...,
\end{equation}
where the hamiltonian $H$ is
\begin{equation}\label{Ham}
  H = \Lambda+G=\frac{1}{2}\sum_{j\geq1}\lambda_j(p_j^2+q_j^2)+\frac{1}{4}\sum_{i,j,k,l}G_{ijkl}q_iq_jq_kq_l,
\end{equation}
\begin{equation}\label{cf}
  G_{ijkl}=\frac{1}{\sqrt{\lambda_i\lambda_j\lambda_k\lambda_l}}\int_{-1}^{1}\phi_i\phi_j\phi_k\phi_l\md x.
\end{equation}

Denoting the invariant $2\times 2$-dimensional linear space by $E$:
\begin{equation*}
  E=\{(u,v)=(q_1\phi_1+q_2\phi_2,p_1\phi_1+p_2\phi_2)\}=\bigcup_{I\in \overline{\mathbb{P}^2}}\mathscr{T}(I),
\end{equation*}
where $\mathbb{P}^2=\{I\in \mathbb{R}^2:I_j>0 ~ \text{for}~ j=1,2\}$ is the positive quadrant in $\mathbb{R}^2$ ,
\begin{equation*}
  \mathscr{T}(I)=\{(u,v):q_j^2+p_j^2=I_j~\text{for}~j=1,2\},
\end{equation*}
then our main theorem is as follows.
\begin{thm}
Consider the nonlinear hamiltonian partial differential equation with boundary condition
 \begin{equation}\label{equation}
\left\{
\begin{aligned}
&u_{tt}-((1-x^2)u_x)_x+mu+ u^3=0,\quad x\in[-1,1]\\
&\lim_{x\rightarrow\pm 1}u(x)(1-x^2)^{\frac{1}{4}}=0.
\end{aligned}
\right.
\end{equation}
If $m\in(0,\frac{1}{4})\cup(\frac{1}{4},\frac{41}{4})$, then there is a set $\mathscr{C}$ in $\mathbb{P}^2$ with positive lebesgue measure, a family of 2-tori
\begin{equation*}
\mathscr{T}[\mathscr{C}]=\bigcup_{I\in \mathscr{C}}\mathscr{T}(I)\subset E
\end{equation*}
over $\mathscr{C}$, as well as a lipschitz continuous embedding into phase space $\mathscr{P}$
\begin{equation*}
\Phi:\mathscr{T}[\mathscr{C}]\hookrightarrow \mathscr{P},
\end{equation*}
which is a higher order perturbation of the inclusion map $\Phi_0:E\hookrightarrow \mathscr{P}$ restricted to $\mathscr{T}[\mathscr{C}]$, such that the restriction of $\Phi$ to each $\mathscr{T}(I)$ in the family is an embedding of a rotational invariant 2-torus for the nonlinear hamiltonian differential equation (\ref{equation}).
\end{thm}

Here are some remarks.
We compare our results with those of P\"{o}schel \cite{poschel1996quasi}. By
and large, the basic idea is the same in reducing the hamiltonian defined by the partial differential equations  to a partial Birkhoff normal form such that the KAM theorem \cite{poschel2006kam} (also see \cite{§Ü§å§Ü§ã§Ú§ßnearly}) is applicable. However, there are several main differences because of the singularity of the differential operator $A$. In  P\"{o}schel \cite{poschel1996quasi} , the differential operator $\hat{A}=-\frac{\md^2}{\md x^2}+m$ with Dirichlet boundary conditions has eigenvalues $\hat{\lambda}_j^2$ and eigenfunction $\hat{\phi}_j$:
\begin{equation*}
  \hat{\lambda}_j^2=j^2+m,\quad\hat{\phi}_j=\sqrt{\frac{2}{\pi}}\sin{jx}.
\end{equation*}
In contrast, the singular differential operator $A$ has, respectively, the eigenvalues and eigenfunctions \begin{equation}\label{eigen}
  \lambda^2_j=2j(2j-1)+m,\quad \phi_j=\sqrt{2j-\frac{1}{2}}P_{2j-1}(x),\quad j=1,2,\ldots
\end{equation}
where $P_j(x)$ are  Legendre polynomials.

On the one hand, under the basis $\{\hat\phi_j\}$ the Hamiltonian of $u^3$ can be written as
\[\hat G(q)=\sum_{ijkl} \hat{G}_{iijj} q_i\,q_j\,q_k\,q_l\]
with
\begin{equation*}
  \hat{G}_{ijkl}=\frac{1}{\sqrt{\hat{\lambda}_i\hat{\lambda}_j\hat{\lambda}_k\hat{\lambda}_l}}\int_0^{\pi}\hat{\phi}_i\hat{\phi}_j\hat{\phi}_k\hat{\phi}_l\md x
\end{equation*}
Since $\hat\phi_j$ is a very simple triangle function $\sqrt{\frac{2}{\pi}}\sin\,x$, it is easy to verify that
\begin{equation}\label{coefficient}
  \hat{G}_{iijj}=\frac{1}{2\pi}\frac{2+\delta_{ij}}{\lambda_i\lambda_j},
\end{equation}
and to fulfill the relationship
\begin{equation}\label{relation}
  \hat{G}_{ijkl}=0
\quad\text{
unless}\quad
  i\pm j\pm k\pm l=0,
\end{equation} where $\delta_{ij}=1 $ when $i=j$ and $\delta_{ij}=0$ when $i\neq j$.
The relationship (\ref{relation}) leads immediately to that the Hamiltonian $\hat G(q)$ is the convolution of $q$ and $q$'s, that is,
\[\hat G(q)=q\star q\star q\star q,\]
form which the regularity of the vector field $X_{\hat G}$ follows.
At the same time, since the coefficients $\hat{G}_{iijj}$ can be explicitly calculated in (\ref{coefficient}), the resonant conditions in both Birkhoff normal form and the KAM theorem can be directly to verified.

 However, on the other hand, under the Legendre basis $\phi_j$'s, the Hamiltonian of $u^3$ can be written as
\[G(q)=\sum_{ijkl} {G}_{iijj} q_i\,q_j\,q_k\,q_l\]
with
\begin{equation*}
  {G}_{ijkl}=\frac{1}{\sqrt{{\lambda}_i{\lambda}_j{\lambda}_k{\lambda}_l}}\int_{-1}^{1}{\phi}_i{\phi}_j{\phi}_k{\phi}_l
  \md x
\end{equation*}
Both the equation (\ref{coefficient}) and the relationship (\ref{relation}) do not hold true any more in this case. Actually, the calculation of the integral $\int_{-1}^{1}{\phi}_i{\phi}_j{\phi}_k{\phi}_l
  \md x$ is not completely solved even in special function theory.  Thus the fulfillment of the regularity of the vector filed $X_G$ and those resonant conditions in both Birkhoff normal form and the KAM theorem are not easy. Section \ref{functionspace} will be devoted to verify the regularity of $X_G$. And the loss of (\ref{coefficient}) accounts for why we choose $m\in(0,\frac{1}{4})\cup(\frac{1}{4},\frac{41}{4})$ and consider only $2$ dimensional KAM tori.

\section{Legendre polynomials and Algebraic Property \label{functionspace}}

In the section, let us introduce some properties about Legendre polynomials $P_n(x)$ first. By using them, we can derive the estimate of $G_{ijkl}$ in next section.

For fixed $n$, the Legendre polynomial $P_n(x)$ is a $n$ order polynomial. It has an usual expression
\begin{equation}\label{def}
  P_n(x)=\sum_{k=0}^{[\frac{n}{2}]}(-1)^k\frac{(2n-2k)!}{2^nk!(n-k)!(n-2k)!}x^{n-2k},
\end{equation}
as well as the Rodrigues's formula
\begin{equation}\label{Def}
  P_n(x)=\frac{1}{2^nn!}\frac{\md^n}{\md x^n}[(x^2-1)^n].
\end{equation}

At the endpoint $x=\pm 1$, it satisfies
\begin{equation}\label{special}
  P_n(1)=1,\quad P_n(-1)=(-1)^n,
\end{equation}
and it has a uniform upper bound
\begin{equation}\label{bound}
  |P_n(x)|\leq 1.
\end{equation}

The recursion formula is important
\begin{equation}\label{recursion}
(n+1)P_{n+1}-(2n+1)xP_n+nP_{n-1}=0.
\end{equation}
A routine computation from (\ref{recursion}) gives rise to
\begin{equation}\label{re1}
  xP_n'-P_{n-1}'=nP_n,
\end{equation}
and
\begin{equation}\label{re2}
  (1-x^2)P_n'=n(P_{n-1}-xP_n).
\end{equation}

From the Rodrigues's formula (\ref{Def}), we get
\begin{equation}\label{pro}
\int_{-1}^{1}x^k P_n(x)\md x =
\left\{\begin{aligned}
&0\quad\quad &k<n\\
&\frac{2^{n+1}(n!)^2}{(2n+1)!}\quad\quad&k=n\\
&0\quad\quad &k>n,k-n\in 2\mathbb{Z}+1\\
&\frac{k!\Gamma(\frac{k}{2}-\frac{n}{2}+\frac{1}{2})}{2^l(k-n)!\Gamma(\frac{k}{2}+\frac{n}{2}+\frac{3}{2})}&k>n,k-n\in 2\mathbb{Z}
\end{aligned}\right.
\end{equation}

A classical formula can express the product of two Legendre polynomials as a sum of such polynomials:
\begin{equation}
P_k(x)P_l(x)=\sum_{m=|k-l|}^{k+l}\frac{\mathbf{A}(s-k)\mathbf{A}(s-l)\mathbf{A}(s-m)}{\mathbf{A}(s)}\frac{2m+1}{2s+1}P_m(x),
\end{equation}
where
\begin{equation*}
 s = \frac{k+l+m}{2}\quad \text{and}\quad \mathbf{A}(n)=\frac{1\cdot3\cdot5\cdot\ldots\cdot(2n-1)}{1\cdot2\cdot3\cdot\ldots\cdot n}=\frac{(2n)!}{2^n(n!)^2}=\frac{1}{2^n}{2n \choose n}.
\end{equation*}
The result can also be expanded in a series using $3j$ symbol as:
\begin{equation}\label{Adam}
P_k(x)P_l(x)=\sum_{m=|k-l|}^{k+l}{k \quad l\quad  m \choose 0\quad 0\quad  0}^2(2m+1)P_m(x),
\end{equation}
where
\begin{equation*}
  {k \quad l\quad  m \choose 0\quad 0\quad  0}^2\triangleq \frac{\mathbf{A}(s-k)\mathbf{A}(s-l)\mathbf{A}(s-m)}{(2s+1)\mathbf{A}(s)}.
\end{equation*}
Thus we could calculate the integral of three Legendre polynomials:
\begin{equation}
\int_{-1}^{1}P_k(x)P_l(x)P_m(x)~\md x=2{k \quad l\quad  m \choose 0\quad 0\quad  0}^2.
\end{equation}

We remark that the result we get in this paper is an extension of the research in special function and refer to \cite{hs1938certain} for details.

Next, let us verify the algebraic property of the function space given below. Employing the result, we can get the regularity of vectorfield in next section.

Let
\begin{equation}\label{sp1}
\mathscr{D}(A^{\frac{s}{2}})=\left\{u\in L^2[-1,1]\Big{|}A^{s/2}u\in L^2[-1,1],\lim_{x\rightarrow \pm 1} u(x)(1-x^2)^{\frac{1}{4}}=0\right\},
\end{equation}
and
\begin{equation}\label{sp2}
   \ell^2_s=\left\{u= \sum_{j\geq 1}u_j \phi_j\in L^2[-1,1]\Big{|}\sum_{j\geq 1}j^{2s}|u_j|^2<\infty\right\},
\end{equation}
where $\|A^{s/2}u\|_{L^2[-1,1]}=\langle A^{s/2}u,A^{s/2}u \rangle^{\frac{1}{2}}=\langle A^su,u \rangle^{\frac{1}{2}}$, and $\langle \cdot,\cdot\rangle$ denotes the usual scalar product in $L^2[-1,1]$. The property
\begin{equation}\label{dd}
  \mathscr{D}(A^{t+1})\subseteq \mathscr{D}(A^t),
\end{equation}
is also necessary.

If the norm of $\ell^2_s$ is defined by $\|u\|_s=(\sum_{j\geq 1}j^{2s}|u_j|^2)^{\frac{1}{2}}$, then the following norms are equivalent
\begin{equation}
  \|A^{\frac{s}{2}}u\|_{L^2[-1,1]}\sim\|u\|_s.
\end{equation}

Remark:On one hand, \\
$\norm{A^{\frac{s}{2}}\sum_{j=1}^nu_j\phi_j}_{L^2[-1,1]}^2\leq \norm{\sum_{j=1}^nu_j\lambda_j^s\phi_j}_{L^2[-1,1]}^2\lesssim\sum_{j=1}^n|u_j|^2j^{2s}\lesssim\norm{u}_s^2\quad n\rightarrow \infty.$\\
On the other hand, \\
$\sum_{j=1}^n|u_j|^2j^{2s}=\int_{-1}^1|\sum_{j=1}^nu_jj^{s}\phi_j|^2\md x\lesssim \int_{-1}^1\left|\sum_{j=1}^n\lambda_j^su_j\phi_j\right|^2\md x\lesssim \int_{-1}^1\left| A^{\frac{s}{2}}\sum_{j=1}^nu_j\phi_j \right|^2\md x\\ \lesssim\norm{A^{\frac{s}{2}}u}_{L^2[-1,1]}^2\quad n\rightarrow \infty.$

 Our theorem is as follows

 \begin{thm}
 If $u(-x)=-u(x), u\in \mathscr{D}(A^2)$,  then
 \begin{equation}\label{a2}
 \norm{A^2 (u^3)}_{L^2[-1,1]}\lesssim \norm{A^2 u}^3_{L^2[-1,1]}.
 \end{equation}
 \end{thm}
\textbf{Proof.}\quad First, we claim that
 \begin{equation}\label{phic}
   |\partial_x\phi_j(x)|\leq \sqrt{2j-1}j(2j-1)\lesssim j^{\frac{5}{2}}.
 \end{equation}
 Let $f(x)=j(1-x^2)-\Big[P_{2j-2}(x)-xP_{2j-1}(x)\Big]$, according to (\ref{re1}), we have
 \begin{equation*}
   f'(x)=2j\Big(P_{2j-1}(x)-x\Big).
 \end{equation*}
 At the critical points $x_0$ such that $P_{2j-1}(x_0)=x_0$, and the endpoint $x=\pm 1$, we find
 \begin{equation*}
 \begin{aligned}
   f(x_0)&=j-P_{2j-2}(x_0)-(j-1)x_0^2\\
   &\geq 0
 \end{aligned}
 \end{equation*}
 and
 \begin{equation*}
   f(-1)=f(1)=0.
 \end{equation*}
 Here we use the property (\ref{special}) and (\ref{bound}). This derives the relationship
 \begin{equation}\label{>}
   j(1-x^2)\geq P_{2j-2}(x)-xP_{2j-1}(x).
 \end{equation}
 Using (\ref{special}) and (\ref{re1}) again, we obtain
 \begin{equation}\label{-1}
   \begin{aligned}
   \lim_{x\rightarrow -1}\frac{P_{2j-2}(x)-xP_{2j-1}(x)}{1-x^2}&=\lim_{x\rightarrow -1}\frac{P_{2j-2}'(x)-P_{2j-1}(x)-xP_{2j-1}'(x)}{-2x}\\
   &=\lim_{x\rightarrow -1}\frac{-2jP_{2j-1}(x)}{-2x}\\
   &=j.
   \end{aligned}
 \end{equation}
 The same method gives
 \begin{equation}\label{1}
   \lim_{x\rightarrow 1}\frac{P_{2j-2}(x)-xP_{2j-1}(x)}{1-x^2}=j.
 \end{equation}
 Combining (\ref{>}),(\ref{-1}),(\ref{1}) and the property (\ref{re2}), we see that
 \begin{equation}
   |P'_{2j-1}(x)|\leq j(2j-1).
 \end{equation}
This leads to (\ref{phic}) because of $\phi_j=\sqrt{2j-\frac{1}{2}}P_{2j-1}$.

 Due to the fact that $u=\sum_{j\geq 1}\frac{q_j(t)}{\sqrt{\lambda_j}}\phi_j\in\mathscr{D}(A^2) $, we have $u\in \ell^2_4$. Using (\ref{phic}), we get
 \begin{equation*}
  \begin{aligned}
   \sum_{j\geq 1}\Big| \frac{q_j(t)}{\sqrt{\lambda_j}}\partial_x\phi_j \Big|&\leq\sum_{j\geq 1}\Big| \frac{q_j(t)}{\sqrt{\lambda_j}}\sqrt{2j-\frac{1}{2}}j(2j-1)\Big|\\
   &\lesssim \Big(\sum_{j\geq 1}\frac{1}{j^3}\Big)^\frac{1}{2}\Big(\sum_{j\geq 1}j^{8}\Big|\frac{q_j(t)}{\sqrt{\lambda_j}}\Big|^2\Big)^{\frac{1}{2}}\\
   &\lesssim\|u\|_{4}\sim \norm{A^2u}_{L^2[-1,1]}.
   \end{aligned}
 \end{equation*}
 Then the sum $\sum_{j\geq 1}\frac{q_j(t)}{\sqrt{\lambda_j}}\partial_x\phi_j(x)$ is absolutely convergent in $\mathbb{R}_+\times [-1,1]$, where $(t,x)\in \mathbb{R}_+\times [-1,1]$. It follows the estimate
 \begin{equation}\label{estimate1}
   \abs{\partial_xu}\lesssim\norm{A^2u}_{L^2[-1,1]},
 \end{equation}
 and
 \begin{equation}\label{estimate2}
   \lim_{x\rightarrow\pm 1}(1-x^2)^{r}\partial_xu=0.\quad(r>0~ \text{in particular} ~ r=1 ~\text{and}~ r=\frac{1}{2})
 \end{equation}

 Since $u(-x)=-u(x)$, we have
\begin{equation}\label{poc}
  \int_{-1}^{1}u\md x=0,\quad\int_{-1}^{1}u^3\md x=0.
\end{equation}
In view of (\ref{estimate1}) , (\ref{estimate2}),(\ref{poc}) and Poincar\'{e} inequality, we get
\begin{equation}\label{pp}
  \norm{u}_{L^2[-1,1]}\lesssim\norm{\partial_xu}_{L^2[-1,1]}\lesssim \norm{A^2u}_{L^2[-1,1]},
\end{equation}

Using Theorem 7.2 (the Gagliardo-Nirenberg Inequality (\ref{gnin})) in Appendix A, for $f(x)=u(x)$,
\begin{equation*}
  k=0,p=\infty,q=2,r=2,m=1,a=\frac{1}{2},
\end{equation*}
we obtain
\begin{equation}\label{ppp}
  \norm{u}_{L^\infty[-1,1]}\lesssim\norm{\partial_xu}_{L^2[-1,1]}^{\frac{1}{2}}\norm{u}_{L^2[-1,1]}^{\frac{1}{2}}\lesssim \norm{A^2u}_{L^2[-1,1]}.
\end{equation}
Then we can conclude
\begin{equation}
  \norm{\mathscr{A}u}_{L^2[-1,1]}\lesssim\norm{A^2u}_{L^2[-1,1]},\quad \norm{\mathscr{A}^2u}_{L^2[-1,1]}\lesssim\norm{A^2u}_{L^2[-1,1]}.
\end{equation}

Noting (\ref{estimate2}) (\ref{poc}) (\ref{pp}) and using Poincar\'{e} inequality, we have the following inequality,
\begin{equation}
\begin{aligned}
  \norm{Au^3}_{L^2[-1,1]}&=\norm{\mathscr{A}u^3+mu^3}_{L^2[-1,1]}\\
  &\leq \norm{\mathscr{A}u^3}_{L^2[-1,1]}+m\norm{u^3}_{L^2[-1,1]}\\
  &\leq \norm{3u^2(1-x^2)\partial_xu}_{H^1[-1,1]}+m\norm{u^3}_{H^1[-1,1]}\\
  &\leq 3\norm{u}_{H^1[-1,1]}^2\norm{(1-x^2)\partial_xu}_{H^1[-1,1]}+m\norm{u}_{{H^1[-1,1]}}^3\\
  &\lesssim\norm{\partial_xu}_{L^2[-1,1]}^2\norm{\partial_x[(1-x^2)\partial_xu]}_{L^2[-1,1]}+m\norm{\partial_xu}_{{L^2[-1,1]}}^3\\
  &\lesssim \norm{A^2u}^2_{L^2[-1,1]}\norm{\mathscr{A}u}_{L^2[-1,1]}+ \norm{A^2u}^3_{L^2[-1,1]}\\
  &\lesssim  \norm{A^2u}^3_{L^2[-1,1]}.
\end{aligned}
\end{equation}
Here, $\norm{f}_{H^1[-1,1]}=\left(\norm{f}_{L^2[-1,1]}^2+\norm{\partial_xf}_{L^2[-1,1]}^2\right)^{\frac{1}{2}}$.

 In 1959, Nirenberg \cite{nirenberg2011elliptic} observed a connection between $L^p$- norms and the H\"{o}lder seminorms $\left[\cdot\right]_\alpha$. Define
\begin{equation}\label{mixnorm}
  \{f\}_\alpha=\left\{
  \begin{aligned}
  &[f]_\alpha=\sup_{x\neq y}\frac{|f(x)-f(y)|}{|x-y|^\alpha}\quad when \quad 0<\alpha<1,\\
  &\norm{f}_{L^p[-1,1]}\quad when \quad \alpha=-\frac{1}{p}\leq 0.
  \end{aligned}\right.
\end{equation}

By  Theorem 7.1 (the General Nirenberg Inequality (\ref{gn1}) ) in Appendix A, for $f(x)=\phi_j(x)$,
\begin{equation*}
\begin{aligned}
  j=0\quad \beta&=\frac{495}{1364}\quad p_\beta=\infty\\
  k=1\quad\alpha&=-\frac{1}{20}\quad p_\alpha=20\\
  \theta=\frac{5885}{9889}\quad\gamma&=-\frac{1}{2}\quad p_\gamma=2
\end{aligned}
\end{equation*}
we have
\begin{equation}
  \left\{\phi_j\right\}_{\frac{495}{1364}}\lesssim\left\{\partial_x\phi_j\right\}_{-\frac{1}{20}}^\frac{5885}{9889}\left\{\phi_j\right\}_{-\frac{1}{2}}^\frac{4004}{9889}.
\end{equation}
Using (\ref{phic}), we obtain
\begin{equation}
\begin{aligned}
  \left[\phi_j\right]_{\frac{495}{1364}}&\lesssim\norm{\partial_x\phi_j}_{L^{20}[-1,1]}^\frac{5885}{9889}\norm{\phi_j}_{L^2[-1,1]}^\frac{4004}{9889}\\
  &\lesssim|\partial_x\phi_j|^\frac{5885}{9889}\\
  &\lesssim j^\frac{29425}{19778}.
\end{aligned}
\end{equation}
It follows that
\begin{equation}\label{aa}
\begin{aligned}
  \left[\mathscr{A}u\right]_{\frac{495}{1364}}&=\left[\mathscr{A}\sum_{j\geq 1}u_j\phi_j\right]_{\frac{495}{1364}}\\
  &\lesssim \sum_{j\geq 1}|u_{j}| j^2\left[\phi_j\right]_{\frac{495}{1364}}\\
  &\lesssim \left(\sum_{j\geq 1}\frac{1}{j^{\frac{20262}{19778}}}\right)^\frac{1}{2}\left(\sum_{j\geq 1}|u_j|^2 j^8\right)^\frac{1}{2}\thicksim \norm{A^2u}_{L^2[-1,1]}
\end{aligned}
\end{equation}
With the help of Theorem 7.1 (the General Nirenberg Inequality (\ref{gn1}) ) in Appendix A again, for $f(x)=(1-x^2)\partial_xu$,
\begin{equation*}
\begin{aligned}
  j=1\quad \beta&=-\frac{1}{4}\quad p_\beta=4\\
  k=1\quad\alpha&=\frac{495}{1364}\quad p_\alpha=\infty\\
  \theta=\frac{31}{42}\quad\gamma&=-\frac{43}{44}\quad p_\gamma=\frac{44}{43}
\end{aligned}
\end{equation*}
we have
\begin{equation}
  \left\{\partial_x[(1-x^2)\partial_xu]\right\}_{-\frac{1}{4}}\lesssim\left\{\partial_x[(1-x^2)\partial_xu]\right\}_{\frac{495}{1364}}^{\frac{31}{42}}\left\{(1-x^2)\partial_xu\right\}_{-\frac{43}{44}}^{\frac{11}{42}}.
\end{equation}
From (\ref{estimate2}) (\ref{aa}), the Sobolev Imbedding theorem $H^1[-1,1]\rightarrow L^{\frac{44}{43}}[-1,1]$ and Poincar\'{e} Inequality, it follows that
\begin{equation}\label{l4}
  \begin{aligned}
  \norm{\mathscr{A}u}_{L^4[-1,1]}&\lesssim\left[\mathscr{A}u\right]_{\frac{495}{1364}}^{\frac{31}{42}}\norm{(1-x^2)\partial_xu}_{L^{\frac{44}{43}}[-1,1]}^{\frac{11}{42}}\\
  &\lesssim \norm{A^2u}_{L^2[-1,1]}^{\frac{31}{42}}\norm{(1-x^2)\partial_xu]}_{H^1[-1,1]}^{\frac{11}{42}}\\
  &\lesssim \norm{A^2u}_{L^2[-1,1]}^{\frac{31}{42}}\norm{\mathscr{A}u}_{L^2[-1,1]}^{\frac{11}{42}}\\
  &\lesssim\norm{A^2u}_{L^2[-1,1]}.
  \end{aligned}
\end{equation}

Then, using (\ref{pp}) (\ref{ppp}) and
\begin{equation}\label{tran}
  (1-x^2)\partial^2_xu=2x\partial_xu-\mathscr{A}u
\end{equation}
we can deduce that
\begin{equation}
\begin{aligned}
  \norm{A^2u^3}_{L^2[-1,1]}&=\norm{\mathscr{A}(Au^3)+mAu^3}_{L^2[-1,1]}\\
  &\lesssim\norm{\mathscr{A}^2u^3}_{L^2[-1,1]}+\norm{A^2u}^3_{L^2[-1,1]}\\
  &\lesssim\norm{(1-x^2)\partial_x^2\left[3 u^2(1-x^2)\partial_xu\right]}_{H^1[-1,1]}+\norm{A^2u}^3_{L^2[-1,1]}\\
  &\lesssim\left\Vert-6\left((1-x^2)\partial_xu\right)\left((1-x^2)(\partial_xu)^2\right)-12xu\left((1-x^2)(\partial_xu)^2\right)\right.\\
  &\left.\quad+18u\left((1-x^2)(\partial_xu)\mathscr{A}u\right)+3u\left(u(1-x^2)\partial_x(\mathscr{A}u)\right)\right\Vert_{H^1[-1,1]}+\norm{A^2u}^3_{L^2[-1,1]}\\
  &\lesssim\norm{(1-x^2)\partial_xu}_{H^1[-1,1]}\norm{(1-x^2)(\partial_xu)^2}_{H^1[-1,1]}+ \norm{A^2u}_{L^2[-1,1]}\norm{(1-x^2)(\partial_xu)^2}_{H^1[-1,1]}\\
  &\quad+\norm{A^2u}_{L^2[-1,1]}\norm{(1-x^2)(\partial_xu)\mathscr{A}u}_{H^1[-1,1]}+\norm{A^2u}_{L^2[-1,1]}\norm{u(1-x^2)\partial_x(\mathscr{A}u)}_{H^1[-1,1]}\\
  &\quad+\norm{A^2u}^3_{L^2[-1,1]}.
\end{aligned}
\end{equation}

 By Poincar\'{e} inequality and (\ref{estimate2}), we obtain
\begin{equation}\label{b1}
  \norm{(1-x^2)\partial_xu}_{H^1[-1,1]}\lesssim\norm{\mathscr{A}u}_{L^2[-1,1]}\lesssim \norm{A^2u}^2_{L^2[-1,1]}.
\end{equation}

Using Poincar\'{e} inequality and (\ref{estimate2}) again with the help of (\ref{estimate1}) and (\ref{tran}) , we have
\begin{equation}\label{b2}
  \begin{aligned}
  &\quad\norm{(1-x^2)(\partial_xu)^2}_{H^1[-1,1]}\\
  &\lesssim\norm{\partial_x\left[\left((1-x^2)(\partial_xu)\right)\partial_xu\right]}_{L^2[-1,1]}\\
  &\lesssim\norm{\partial_xu\partial_x\left((1-x^2)(\partial_xu)\right)}_{L^2[-1,1]}+\norm{(\partial_xu)(1-x^2)\partial_x^2u}_{L^2[-1,1]}\\
  &\lesssim\norm{A^2u}_{L^2[-1,1]}\norm{\mathscr{A}u}_{L^2[-1,1]}+\norm{\partial_xu\left(2x\partial_xu-\mathscr{A}u\right)}_{L^2[-1,1]}\\
  &\lesssim\norm{A^2u}_{L^2[-1,1]}\norm{Au}_{L^2[-1,1]}+\norm{A^2u}^2_{L^2[-1,1]}\\
  &\lesssim\norm{A^2u}^2_{L^2[-1,1]}.
  \end{aligned}
\end{equation}

From (\ref{dd}), we know for $u\in \mathscr{D}(A^2), Au\in \mathscr{D}(A)$, then $\lim_{x\rightarrow \pm1}(1-x^2)^{\frac{1}{4}}Au=0$. So we have
\begin{equation}
  \lim_{x\rightarrow \pm1}(1-x^2)\mathscr{A}u=0
\end{equation}
Using (\ref{estimate2}), Cauchy Inequality and Poincar\'{e } Inequality, we get
\begin{equation}\label{ave}
  \begin{aligned}
  \left|\int_{-1}^1(1-x^2)\partial_x(\mathscr{A}u)\md x\right|&=\left|\int_{-1}^1 \partial_x[(1-x^2)\mathscr{A}u] \md x-\int_{-1}^1 [\partial_x(1-x^2)]\mathscr{A}u\md x\right|\\
  &=2\left|\int_{-1}^1x\left[2x\partial_xu-(1-x^2)\partial_x^2u\right]\md x\right|\\
  &=2\left|\int_{-1}^1 2x^2\partial_xu\md x-\int_{-1}^1x(1-x^2)\partial_x^2u\md x\right|\\
  &=2\left|\int_{-1}^1 2x^2\partial_xu\md x-\int_{-1}^1\partial_x\left(x(1-x^2)\partial_xu\right)\md x+\int_{-1}^1\partial_x\left(x(1-x^2)\right)\partial_xu\md x\right|\\
  &=2\left|\int_{-1}^1 (1-x^2)\partial_xu \md x\right|\\
  &\lesssim \norm{(1-x^2)\partial_xu}_{L^2[-1,1]}\\
  &\lesssim\norm{\mathscr{A}u}_{L^2[-1,1]}.
  \end{aligned}
\end{equation}
Then it follows from (\ref{ave})
\begin{equation}\label{pave}
  \begin{aligned}
  \norm{(1-x^2)\partial_x(\mathscr{A}u)}_{L^2[-1,1]}&\lesssim \norm{\partial_x\left((1-x^2)\partial_x(\mathscr{A}u)\right)}_{L^2[-1,1]}+\left|\int_{-1}^1(1-x^2)\partial_x(\mathscr{A}u)\md x\right|\\
  &\lesssim \norm{\mathscr{A}^2u}_{L^2[-1,1]}+\norm{\mathscr{A}u}_{L^2[-1,1]}\\
  &\lesssim \norm{A^2u}_{L^2[-1,1]}
  \end{aligned}
\end{equation}

By Poincar\'{e} Inequality with the help of (\ref{estimate1}) (\ref{estimate2}) (\ref{l4}) and (\ref{pave}), we have
\begin{equation}\label{b3}
  \begin{aligned}
 &\quad\norm{(1-x^2)(\partial_xu)\mathscr{A}u}_{H^1[-1,1]}\\
 &\lesssim\norm{\partial_x\left((1-x^2)(\partial_xu)\mathscr{A}u\right)}_{L^2[-1,1]}\\
  &\lesssim\norm{\left(\partial_x(1-x^2)(\partial_xu)\right)\mathscr{A}u}_{L^2[-1,1]}+\norm{(\partial_xu)(1-x^2)\partial_x(\mathscr{A}u)}_{L^2[-1,1]}\\
  &\lesssim\norm{(\mathscr{A}u)^2}_{L^2[-1,1]}+\norm{A^2u}_{L^2[-1,1]}\norm{(1-x^2)\partial_x(\mathscr{A}u)}_{L^2[-1,1]}\\
  &\lesssim\norm{\mathscr{A}u}_{L^4[-1,1]}^2+\norm{A^2u}_{L^2[-1,1]}^2\\
  &\lesssim \norm{A^2u}_{L^2[-1,1]}^2.
  \end{aligned}
\end{equation}

Then we obtain the following inequality by (\ref{estimate1})  (\ref{ppp}) (\ref{pave}) and Poincar\'{e} Inequality
\begin{equation}\label{b4}
  \begin{aligned}
  \norm{u(1-x^2)\partial_x(\mathscr{A}u)}_{H^1[-1,1]}&\lesssim \norm{\partial_x\left(u(1-x^2)\partial_x(\mathscr{A}u)\right)}_{L^2[-1,1]}\\
  &\lesssim \norm{u\partial_x\left((1-x^2)\partial_x(\mathscr{A}u)\right)}_{L^2[-1,1]}+\norm{(\partial_xu)(1-x^2)\partial_x(\mathscr{A}u)}_{L^2[-1,1]}\\
  &\lesssim \norm{u}_{L^\infty[-1,1]}\norm{\mathscr{A}^2u}_{L^2[-1,1]}+\norm{A^2u}_{L^2[-1,1]}^2\\
  &\lesssim \norm{A^2u}_{L^2[-1,1]}^2
  \end{aligned}
\end{equation}

In view of (\ref{b1})(\ref{b2})(\ref{b3})(\ref{b4}), we proof the inequality (\ref{a2}).
\qed
 \section{The hamiltonian }
  \emph{3.1 the regularity of vectorfield.}\\

 From introduction, we have already obtain the hamiltonian (\ref{Ham})
\begin{equation*}
H = \Lambda+G=\frac{1}{2}\sum_{j\geq1}\lambda_j(p_j^2+q_j^2)+\frac{1}{4}\sum_{i,j,k,l}G_{ijkl}q_iq_jq_kq_l,
\end{equation*}
with equations of motions(\ref{qp})
\begin{equation*}
\dot{q_j}=\frac{\partial H}{\partial p_j}=\lambda_jp_j,\quad \dot{p_j}=\frac{\partial H}{\partial q_j}=-\lambda_jq_j-\frac{\partial G}{\partial q_j},
\end{equation*}
in some neighbourhood of the origin in the Hilbert space $\ell^2_s\times \ell^2_s$ with standard symplectic structure $\sum_j\md q_j \wedge \md p_j$. Then we have the following lemma
\begin{lem}
For $u (=\sum_{j\geq 1} \frac{q_j}{\sqrt{\lambda_j}}\phi_j)\in \mathscr{D}(A^2)$ or $q=(q_j)_{j\geq 1}\in\ell^2_\frac{7}{2}$, the gradient $G_q=(\frac{\partial G}{\partial q_j})_{j\geq1}$ is real analytic as a map from some neighbourhood of the origin in $\ell^2_\frac{7}{2}$ into $\ell^2_\frac{9}{2}$, with
\begin{equation}
\norm{G_q}_\frac{9}{2}= O(\norm{q}^3_\frac{7}{2}).
\end{equation}
\end{lem}

\textbf{ Proof.} From the notation
\begin{equation}
\frac{\partial G}{\partial q_j}=\frac{1}{\sqrt{\lambda_j}}\langle u^3,\phi_j\rangle,
\end{equation}
we have
\begin{equation}
\begin{aligned}
\norm{G_q}_\frac{9}{2}&\lesssim \norm{(\langle u^3,\phi_j \rangle)_{j\geq 1}}_4\\
&\sim\norm{A^2u^3}_{L^2[-1,1]}\\
&\lesssim \norm{A^2u}_{L^2[-1,1]}^3\\
&\lesssim \norm{(\langle u,\phi_j \rangle)_{j\geq 1}}_4^3\\
&\lesssim  \norm{q}_\frac{7}{2}^3.
\end{aligned}
\end{equation}

Since $G$ is independent of $p$, the associated hamiltonian vectorfield,
\begin{equation}
X_G=\sum_{j\geq 1}\left(\frac{\partial G}{\partial p_j}\frac{\partial}{\partial q_j}-\frac{\partial G}{\partial q_j}\frac{\partial}{\partial p_j}\right),
\end{equation}
is smoothing of order 1. By contrast, $X_\Lambda$ is unbounded of order 1.\\\qed

\emph{3.2 The Legendre sequences}\\

 It is necessary to make clear the coefficient $G_{ijkl} (\ref{cf})$ in hamiltonian $H$. In particular $G_{iijj}$ . Then we acquire the property of Legendre sequences denoted by $\mathbf{P}(m,n)=\int_{-1}^{1}P_mP_mP_nP_n\md x, m,n \in \mathbb{N}$ below.
\begin{thm}\textbf{(Legendre sequences)}\quad The Legendre sequences $\mathbf{P}(m,n)$ satisfy the following recursion formula
\begin{equation}
\begin{aligned}
\mathbf{P}(m+1,n)=&\alpha_m^{n-1}\mathbf{P}(m,n-1)+\alpha_m^n\mathbf{P}(m,n)+\alpha_m^{n+1}\mathbf{P}(m,n+1)\\
-&\alpha_{m-1}^{n-1}\mathbf{P}(m-1,n-1)-\alpha_{m-1}^n\mathbf{P}(m-1,n)-\alpha_{m-1}^{n+1}\mathbf{P}(m-1,n+1) \\
+&\alpha_{m-2}^n\mathbf{P}(m-2,n),
\end{aligned}
\end{equation}
where
\begin{equation*}
\alpha_m^{n-1}=\left[\frac{(2m+1)n}{(m+1)(2n+1)}\right]^2,
\end{equation*}
\begin{equation*}
\alpha_m^n=\frac{m(2m+1)(n^2+n-m^2-m)}{(m+1)^3(2m-1)}\left[\frac{(m-1)m}{n^2+n-m^2+m}+\frac{2}{(2n+1)^2} \right],
\end{equation*}
\begin{equation*}
\alpha_m^{n+1}=\left[\frac{(2m+1)(n+1)}{(m+1)(2n+1)}\right]^2,
\end{equation*}
\begin{equation*}
\alpha_{m-1}^{n-1}=\frac{(m-1)(2m-1)(2m+1)(n^2+n-m^2-m)n^2}{(m+1)^3(n^2+n-m^2+m)(2n+1)^2},
\end{equation*}
\begin{equation*}
\alpha_{m-1}^n=\frac{2m(n^2+n-m^2-m)}{(m+1)^3(2n+1)^2}+\left(\frac{m}{m+1}\right)^2,
\end{equation*}
\begin{equation*}
\alpha_{m-1}^{n+1}=\frac{(m-1)(2m-1)(2m+1)(n^2+n-m^2-m)(n+1)^2}{(m+1)^3(n^2+n-m^2+m)(2n+1)^2},
\end{equation*}
\begin{equation*}
\alpha_{m-2}^{n}=\frac{(m-1)^3(2m+1)(n^2+n-m^2-m)}{(m+1)^3(2m-1)(n^2+n-m^2+m)},
\end{equation*}
satisfying
\begin{equation}\label{coe}
  \alpha_m^{n-1}+\alpha_m^n+\alpha_m^{n+1}-\alpha_{m-1}^{n-1}-\alpha_{m-1}^n-\alpha_{m-1}^{n+1}+\alpha_{m-2}^n\equiv1,
\end{equation}
with
\begin{equation}\label{PN0}
\mathbf{P}(0,n)=\frac{2}{2n+1},
\end{equation}
\begin{equation}\label{PN1}
\mathbf{P}(1,n)=\frac{2(2n^2+2n-1)}{(2n-1)(2n+1)(2n+3)},
\end{equation}
\begin{equation}\label{PN2}
\mathbf{P}(2,n)=\frac{11n^4+22n^3-31n^2-42n+18}{(2n-3)(2n-1)(2n+1)(2n+3)(2n+5)},
\end{equation}
\begin{equation}\label{PN3}
  \mathbf{P}(3,n)=\frac{34n^6+102n^5-305n^4-780n^3+703n^2+1110n-450}{(2n-5)(2n-3)(2n-1)(2n+1)(2n+3)(2n+5)(2n+7)}.
\end{equation}

Moreover, if $\int_{-1}^{1}P_i(x)P_j(x)P_k(x)P_l(x)\md x\not=0$, then we obtain the estimate of the following integral,
\begin{equation}\label{integral}
\int_{-1}^{1}P_i(x)P_j(x)P_k(x)P_l(x)\md x\lesssim \frac{1}{\sqrt{i+\frac{1}{2}}\sqrt{j+\frac{1}{2}}\sqrt{k+\frac{1}{2}}\sqrt{l+\frac{1}{2}}}
 \end{equation}
or in short
\begin{equation}
\mathbf{P}(m,n)\lesssim \frac{1}{mn}.
\end{equation}

 In particular, there exist an absolute constant $C>0$ such that
 \begin{equation}\label{basis}
0<\int_{-1}^{1}\phi_i\phi_j\phi_k\phi_l \md x\leq C
 \end{equation}
\end{thm}
The proof is left in section 6.

Using the property of the Legendre polynomials (\ref{pro}), we can obtain the property about $G_{ijkl}$ which we need in the next section.
\begin{lem}
Assume $0<i\leq j\leq k\leq l$ and $i+j+k\geq l$, the coefficients $G_{ijkl}=0$ unless $i\pm j\pm k\pm l\in 2\mathbb{Z}$.
\end{lem}

\textbf{Proof.} From the definition of $G_{ijkl}$(\ref{cf}) and $\phi_i$(\ref{eigen})(\ref{def}), we know the product of $\phi_i\phi_j\phi_k$ is a polynomial like $f(x)= \alpha_1 x^{i+j+k}+\alpha_2 x^{i+j+k-2}+\ldots$. Then due to the property (\ref{pro}), under the assumption $0<i\leq j\leq k\leq l$ and $i+j+k\geq l$,we can check
\begin{equation*}
 \begin{aligned}
  G_{ijkl}&=\frac{1}{\sqrt{\lambda_i\lambda_j\lambda_k\lambda_l}}\int_{-1}^{1}f(x)\phi_l\md x\\
  &=\frac{\alpha_1 \sqrt{2l-\frac{1}{2}}}{\sqrt{\lambda_i\lambda_j\lambda_k\lambda_l}}\int_{-1}^1x^{i+j+k}P_l\md x+\frac{\alpha_2  \sqrt{2l-\frac{1}{2}}}{\sqrt{\lambda_i\lambda_j\lambda_k\lambda_l}}\int_{-1}^1x^{i+j+k-2}P_l\md x+\ldots\\
  &=0
 \end{aligned}
\end{equation*}
unless
\begin{equation*}
  i\pm j\pm k\pm l\in 2\mathbb{Z}.
\end{equation*}
 \qed

 \section{Partial Birkhoff normal form }
Next we introduce complex coordinates
\begin{equation}
z_j=\frac{1}{\sqrt{2}}(q_j+\mi p_j),\quad \overline{z}_j=\frac{1}{\sqrt{2}}(q_j-\mi p_j).
\end{equation}
Then we obtain a real analytic hamiltonian $H=\sum_j\lambda_j|z_j|^2+\ldots$ on the complex Hilbert space $\ell^2_s$ with symplectic structure $\mi \sum_j \md z_j \wedge \md \overline{z}_j$.

In the following, $A(\ell^2_s,\ell^2_{s+1})$ denotes the class of all real analytic maps from some neighbourhood of the origin in $\ell^2_s$ into $\ell^2_{s+1}$. Thus we can also obtain the main proposition like P\"{o}schel \cite{poschel1996quasi}  but the handling of the small denominator is more complex.\\

\begin{prop}
For the 2-dimension of the invariant tori and $m\in(0,\frac{1}{4})\cup (\frac{1}{4},\frac{41}{4})$, there exists a real analytic, symplectic change of coordinates $\Gamma$ in some neighbourhood of the origin in $\ell^2_\frac{7}{2}$ that takes the hamiltonian $H=\Lambda+G$ with nonlinearity into
\begin{equation*}
H\circ \Gamma =\Lambda+\bar{G}+\hat{G}+K,
\end{equation*}
where $X_{\bar{G}},X_{\hat{G}},X_K\in A(\ell^2_\frac{7}{2},\ell^2_{\frac{9}{2}})$,
\begin{equation*}
\bar{G}=\frac{1}{2}\sum_{\min(i,j)\leq2}\overline{G}_{ij}|z_i|^2|z_j|^2
\end{equation*}
with uniquely determined coefficient, and
\begin{equation*}
|\hat{G}|=O(\norm{\hat{z}}^4_\frac{7}{2}),\quad|K|=O(\norm{z}^6_\frac{7}{2}),\quad\hat{z}=(z_{3},z_{4},\ldots).
\end{equation*}
Moreover, the neighbourhood can be chosen uniformly for every compact $m-$interval in $(0,\frac{1}{4})\cup (\frac{1}{4},\frac{41}{4})$, and the dependence of ~$\Gamma$ on $m$ is real analytic.
\end{prop}

\begin{rem}
Here $m=\frac{1}{4}$ and $m=\frac{41}{4}$ play an important pole in controlling the small divisor. For $m=\frac{1}{4}$, the resonance occurs when we check if $\lambda_i+\lambda_j+\lambda_k+\lambda_l=0$. For $m=\frac{41}{4}$, one can find the answer in (\ref{m}) of Lemma 4.3.
\end{rem}

Thus, the hamiltonian $\Lambda+\overline{G}$ is integrable with integrals $|z_j|^2,j=1,2$, while the not-normalized fourth order term $\hat{G}$ is not integrable, but independent of the first 2 modes.\\

\textbf{Proof of property.}\quad  Let us introduce another set of coordinates $(\ldots,w_{-2},w_{-1},w_1,w_2,\ldots)$ in $\ell^2_s$ by setting $z_j=w_j,\bar{z}_j=w_{-j}$. The hamiltonian under consideration then reads
\begin{equation*}
\begin{aligned}
H=&\Lambda + G\\
=&\sum_{j\geq1}\lambda_jz_j\bar{z}_j+\frac{1}{4}\sum_{i,j,k,l}G_{ijkl}(z_i+\bar{z}_i)\ldots(z_l+\bar{z}_l)\\
=&\sum_{j\geq1}\lambda_jw_jw_{-j}+{\sum_{i,j,k,l}}^{\prime}G_{ijkl}w_iw_jw_kw_l.
\end{aligned}
\end{equation*}
The prime indicates that the subscripted indices run through all nonzero integers. The coefficients are defined for arbitrary integers by setting $G_{ijkl}=G_{|i||j||k||l|}$.

Formally, the transformation $\Gamma$ is obtained as the time-1-map of the flow of a hamiltonian vectorfield $X_F$ given by a hamiltonian
\begin{equation*}
F={\sum_{i,j,k,l}}^{\prime}F_{ijkl}w_iw_jw_kw_l,
\end{equation*}
with coefficients
\begin{equation}\label{small}
\mi F_{ijkl}=\left\{
\begin{aligned}
&\frac{G_{ijkl}}{\lambda^{\prime}_i+\lambda^{\prime}_j+\lambda^{\prime}_k+\lambda^{\prime}_l} \quad \text{for}& (i,j,k,l)\in \mathscr{L}\backslash \mathscr{N},\\
&0 &\text{otherwise}.
\end{aligned}\right.
\end{equation}
Here, $\lambda_j^{\prime}=$sgn$j\cdot\lambda_{|j|}$,
\begin{equation*}
\mathscr{L}=\left\{(i,j,k,l)\in\mathbb{Z}^4:0\neq \min(|i|,\ldots,|l|)\leq 2 \right\},
\end{equation*}
and $\mathscr{N}\subset\mathscr{L}$ is the subset of all $(i,j,k,l)\equiv(p,-p,q,-q)$. That is, they are of the form $(p,-p,q,-q)$ or some permutation of it.

Next, we will estimate the denominator $\lambda^{\prime}_i+\lambda^{\prime}_j+\lambda^{\prime}_k+\lambda^{\prime}_l$ to ensure the correction of the definition of (\ref{small}), the proof of the lemma is left at the end of this section.
\begin{lem}
If $i,j,k,l$ are non-zero integers, such that $i\pm j\pm k\pm l\in 2\mathbb{Z}$, but $(i,j,k,l)\not\equiv(p,-p,q,-q)$, then
\begin{equation*}
|\lambda^{\prime}_i+\lambda^{\prime}_j+\lambda^{\prime}_k+\lambda^{\prime}_l|\gtrsim \sigma(m,n)=\sigma>0, \quad n=\min(|i|,\ldots,|l|),
\end{equation*}
Hence, the denominators in (\ref{small}) are uniformly bounded away from zero on every compact $m-$interval in $(0,\frac{1}{4})\cup(\frac{1}{4},\frac{41}{4})$. The notation $\sigma(m,n)$ is defined as
\begin{equation*}
\sigma(m,n)=\min\left\{W(m,n),V(m,n)\right\},
\end{equation*}
where
\begin{equation*}
\begin{aligned}
W(m,n)=&\\
=&\min\left\{
\begin{aligned}
&\frac{(1-4m)(n+m)}{2(\sqrt{n(n+1)+m}+n)(\sqrt{n(n+1)+m}+n+2m)}\quad&m\in(0,\frac{1}{4})\\
&2-\frac{4m-1}{4\sqrt{n(n+1)+m}+4n+2}\quad&m\in(\frac{1}{4},\frac{41}{4})
\end{aligned}\right\},
\end{aligned}
\end{equation*}
and
\begin{equation*}
V(m,n)=\min_{m\in(0,\frac{1}{4})\cup(\frac{1}{4},\frac{41}{4})}\left\{\frac{m}{\sqrt{n(n+1)+m}},\frac{n}{\sqrt{m+2}},{\frac{4m-1}{4(n(n+1)+m)^{\frac{3}{2}}}} \right\}.
\end{equation*}
\end{lem}

We continue the proof of the property. Expanding at $t=0$ and using Taylor's formula we formally obtain
\begin{equation*}
\begin{aligned}
H\circ \Gamma=&H\circ X^t_F|_{t=1}\\
=&H+\left\{H,F\right\}+\int_{0}^{1}(1-t)\left\{\left\{H,F \right\},F \right\}\circ X^t_F\md t\\
=&\Lambda+G+\left\{\Lambda,F \right\}+\left\{G,F \right\}+\int_0^1(1-t)\left\{\left\{H,F \right\},F \right\}\circ X^t_F\md t,
\end{aligned}
\end{equation*}
where $\left\{H,F \right\}$ denotes the Poisson bracket of $H$ and $F$. The last line consists of terms of order six or more in $w$ and constitutes the higher order term $K$. In the second to last line,
\begin{equation*}
\left\{\Lambda,F \right\}=-\mi{\sum_{i,j,k,l}}^{\prime}(\lambda^{\prime}_i+\lambda^{\prime}_j+\lambda^{\prime}_k+\lambda^{\prime}_l)F_{ijkl}w_iw_jw_kw_l,
\end{equation*}
hence
\begin{equation*}
G+\left\{\Lambda,F\right\}=\sum_{(i,j,k,l)\in\mathscr{N}}+\sum_{(i,j,k,l)\not\in\mathscr{L}}G_{ijkl}w_iw_jw_kw_l=\bar{G}+\hat{G}.
\end{equation*}
Re-introducing the notations $z_j,\bar{z}_j$ and counting multiplicities, we obtain that
\begin{equation}\label{G1}
\bar{G}=\frac{1}{2}\sum_{\min{(i,j)}\leq2}\bar{G}_{ij}|z_i|^2|z_j|^2,
\end{equation}
with
\begin{equation}\label{G2}
\bar{G}_{ij}=\left\{
\begin{aligned}
&\frac{2(4i-1)(4j-1)}{\lambda_i\lambda_j}\mathbf{P}(2i-1,2j-1)\quad i\neq j,\\
&\frac{(4i-1)(4j-1)}{\lambda_i\lambda_j}\mathbf{P}(2i-1,2j-1)\quad i= j.
\end{aligned}\right.
\end{equation}
Thus, we have $H\circ \Gamma = \Lambda+\bar{G}+\hat{G}+K$ as claimed formally.

To prove analyticity and regularity of the preceding transformation, we first show that
\begin{equation*}
X_F\in A(\ell^2_\frac{7}{2},\ell^2_{\frac{9}{2}}).
\end{equation*}
Assume a `` \emph{threshold function}''
\begin{equation}\label{threshold}
  \tilde{F}=\int_{-1}^1v^4\md x,
\end{equation}
where
\begin{equation*}
  v=\frac{1}{\sigma^{\frac{1}{4}}}\sum_{j}\tilde{w}_j\phi_j\quad \text{and}\quad \tilde{w}_j=\frac{|w_j|+|w_{-j}|}{\sqrt{|j|}}.
\end{equation*}
The natation $\sigma$ comes from the estimate of denominator $|\lambda^{\prime}_i+\lambda^{\prime}_j+\lambda^{\prime}_k+\lambda^{\prime}_l|\gtrsim \sigma$.

It is easy to check that, by (\ref{basis}), the integral of $\int_{-1}^{1}\phi_i\phi_j\phi_k\phi_l \md x$ is uniform bounded.
\begin{equation}\label{ine}
\begin{aligned}
\left|\frac{\partial F}{\partial w_l}\right|& \leq{\sum_{i,j,k}}^{\prime}|F_{ijkl}||w_iw_jw_k|\\
& \lesssim \frac{1}{\sigma\sqrt{|l|}}{\sum_{i,j,k}}^{\prime}\left|\int_{-1}^{1}\phi_i\phi_j\phi_k\phi_l \md x\right|\frac{|w_iw_jw_k|}{\sqrt{|ijk|}}\\
&\lesssim \frac{1}{\sigma\sqrt{|l|}}{\sum_{i,j,k}}^{\prime}\frac{|w_iw_jw_k|}{\sqrt{|ijk|}}
\end{aligned}
\end{equation}
while $\tilde{F}=\frac{1}{\sigma}\sum_{i,j,k,l}\tilde{w}_i\tilde{w}_j\tilde{w}_k\tilde{w}_l\int_{-1}^{1}\phi_i\phi_j\phi_k\phi_l\md x$.

Then it follows that
\begin{equation*}
\begin{aligned}
\left|\frac{\partial \tilde{F}}{\partial w_l}\right|&=\left|\frac{\partial \tilde{F}}{\partial \tilde{w}_l}\cdot\frac{\partial\tilde{w}_l}{\partial w_l}\right|\\
&=\left| \frac{1}{\sigma}\sum_{i,j,k}\tilde{w}_i\tilde{w}_j\tilde{w}_k\int_{-1}^{1}\phi_i\phi_j\phi_k\phi_l \md x \cdot \frac{1}{\sqrt{|l|}}\frac{\bar{w_{l}}}{|w_{l}|}\right|\\
&\gtrsim \frac{1}{\sigma\sqrt{|l|}}\sum_{i,j,k}\left|\int_{-1}^{1}\phi_i\phi_j\phi_k\phi_l \md x\right|\frac{(|w_i|+|w_{-i}|)(|w_j|+|w_{-j}|)(|w_k|+|w_{-k}|)}{\sqrt{|ijk|}}\\
&\gtrsim \frac{1}{\sigma\sqrt{|l|}}{\sum_{i,j,k}}^{\prime}\left|\int_{-1}^{1}\phi_i\phi_j\phi_k\phi_l \md x\right|\frac{|w_iw_jw_k|}{\sqrt{|ijk|}}.
\end{aligned}
\end{equation*}
Hence, the second inequality of (\ref{ine}) implies that
\begin{equation*}
\left|\frac{\partial F}{\partial w_l}\right|\lesssim \left|\frac{\partial \tilde{F}}{\partial w_l}\right|,
\end{equation*}
which means
\begin{equation*}
\norm{F_w}_{\frac{9}{2}}\lesssim \norm{\tilde{F}_w}_{\frac{9}{2}}.
\end{equation*}

On the other hand,
\begin{equation*}
\begin{aligned}
 \left|\frac{\partial \tilde{F}}{\partial w_l}\right|&=\left| \int_{-1}^{1}4v^3\frac{\partial v}{\partial \tilde{w}_l}\cdot\frac{\partial \tilde{w}_l}{\partial w_l}\md x \right|\\
 &=\left|\int_{-1}^{1}4v^3\frac{1}{\sigma^{\frac{1}{4}}}\phi_l\cdot\frac{1}{\sqrt{|l|}}\frac{\bar{w_{l}}}{|w_{l}|}\md x\right|\\
 &\lesssim \frac{1}{\sqrt{|l|}}|\langle v^3, \phi_l \rangle|,
\end{aligned}
\end{equation*}
hence,
\begin{equation*}
\begin{aligned}
\norm{\tilde{F}_w}_\frac{9}{2}=&\left(\sum_ll^{9}\left|\frac{1}{\sqrt{l}}\langle v^3,\phi_l \rangle\right|^2\right)^{\frac{1}{2}}\\
=&\left(\sum_{l}l^{4}|\langle v^3,\phi_l \rangle|^2\right)^{\frac{1}{2}}\\
=&\norm{v^3}_{4}\thicksim \norm{A^{2}v^3}_{L^2[-1,1]}\\
\leq& \norm{A^{2}v}_{L^2[-1,1]}^3\\
\leq&  \norm{v}_{4}^3\\
\thicksim&\left(\sum_l{l^{8}\left|\frac{|w_l|}{\sqrt{|l|}} \right|^{2}}\right)^\frac{3}{2}\thicksim\norm{w}_\frac{7}{2}^3.
\end{aligned}
\end{equation*}
In the end, we obtain
\begin{equation*}
\norm{F_w}_{\frac{9}{2}}\lesssim\norm{w}_\frac{7}{2}^3.
\end{equation*}
The analyticity of $F_w$ follows from the analyticity of each component function and its local boundedness.

In a sufficiently small neighbourhood of the origin in $\ell^2_\frac{7}{2}$, the time-1-map $X_F^t|_{t=1}$ is well defined and gives rise to a real analytic symplectic change of coordinates $\Gamma$ with the estimates
\begin{equation*}
\norm{\Gamma-id}_{\frac{9}{2}}=O(\norm{w}_\frac{7}{2}^3),\quad \norm{D\Gamma-I}_{\frac{9}{2},\frac{7}{2}}^{\mathrm{op}}=O(\norm{w}_\frac{7}{2}^2),
\end{equation*}
where the operator norm $\norm{\cdot}_{r,s}^{\mathrm{op}}$ is defined by
\begin{equation*}
\norm{A}_{r,s}^{\mathrm{op}}=\sup_{w\neq 0}\frac{\norm{Aw}_r}{\norm{w}_s}.
\end{equation*}
Obviously,$\norm{D\Gamma-I}_{\frac{9}{2},\frac{9}{2}}^{\mathrm{op}}\leq\norm{D\Gamma-I}_{\frac{9}{2},\frac{7}{2}}^{\mathrm{op}}$, whence in a sufficiently small neighbourhood of the origin, $D\Gamma$ defines an isomorphism of $\ell^2_{s+1}$. It follows that with $X_H\in A(\ell^2_\frac{7}{2},\ell^2_{\frac{9}{2}})$,
\begin{equation*}
  \Gamma^{*}X_H=D\Gamma^{-1}X_H\circ\Gamma=X_{H\circ\Gamma}\in A(\ell^2_\frac{7}{2},\ell^2_{\frac{9}{2}}).
\end{equation*}
The same holds for the Lie bracket: the boundedness of $\norm{DX_F}_{\frac{9}{2},\frac{7}{2}}^{\mathrm{op}}$ implies that
\begin{equation*}
[X_F,X_H]=X_{\{H,F\}}\in A(\ell^2_\frac{7}{2},\ell^2_{\frac{9}{2}}).
\end{equation*}
These two facts show that $X_K\in A(\ell^2_{\frac{7}{2}},\ell^2_{\frac{9}{2}})$. The analogue claims for $X_{\bar{G}}$ and $X_{\hat{G}}$ are obvious.
\qed\\

\textbf{Proof of Lemma 4.2 }\quad In fact, we want to prove there exists the lower bound of $\lambda^{\prime}_i+\lambda^{\prime}_j+\lambda^{\prime}_k+\lambda^{\prime}_l$, it does not matter to use the renumbered notation $\lambda_j^2=j(j+1)+m$ instead of $\lambda^2_j=2j(2j-1)+m$. This also makes it easier to use \textbf{Lemma 4 in P\"{o}schel } \cite{poschel1996quasi}  in our proof.
\\

\textbf{Case1.} Assume $0<m<\frac{1}{4}$,
using a convenient mark $\sigma_h=\mathrm{sgn} h$, we can write $\delta=\lambda^{\prime}_i+\lambda^{\prime}_j+\lambda^{\prime}_k+\lambda^{\prime}_l=\sigma_i\lambda_{|i|}+\sigma_j\lambda_{|j|}+\sigma_k\lambda_{|k|}+\sigma_l\lambda_{|l|}$ . The skill in \cite{gao2012quasi} will be used.

If $\sigma_i|i|+\sigma_j|j|+\sigma_k|k|+\sigma_l|l|\neq 0$, then $|\sigma_i|i|+\sigma_j|j|+\sigma_k|k|+\sigma_l|l||\geq 2$. Since
\begin{equation*}
\begin{aligned}
\lambda_{|h|}=&|h|+\sqrt{|h|(|h|+1)+m}-|h|\\
=&|h|+\frac{|h|+m}{\sqrt{|h|(|h|+1)+m}+|h|}
\end{aligned}
\end{equation*}
and $f(x)=\frac{x+m}{\sqrt{x(x+1)+m}+x}<\frac{1}{2}$, $(x\geq 0, 0<m<\frac{1}{4})$,
we have
\begin{equation*}
\begin{aligned}
|\delta|=&|\lambda^{\prime}_i+\lambda^{\prime}_j+\lambda^{\prime}_k+\lambda^{\prime}_l|\\
=&|\sigma_i\lambda_{|i|}+\sigma_j\lambda_{|j|}+\sigma_k\lambda_{|k|}+\sigma_l\lambda_{|l|}|\\
\geq&|\sigma_i|i|+\sigma_j|j|+\sigma_k|k|+\sigma_l|l||\\
-&\left( \frac{|i|+m}{\sqrt{|i|(|i|+1)+m}+|i|} +\frac{|j|+m}{\sqrt{|j|(|j|+1)+m}+|j|} +\frac{|k|+m}{\sqrt{|k|(|k|+1)+m}+|k|} +\frac{|l|+m}{\sqrt{|l|(|l|+1)+m}+|l|} \right)\\
\geq&\frac{1}{2}-\frac{|i|+m}{\sqrt{|i|(|i|+1)+m}+|i|}\\
=&\frac{(1-4m)(|i|+m)}{2(\sqrt{|i|(|i|+1)}+|i|)(\sqrt{|i|(|i|+1)+m}+|i|+2m)}.
\end{aligned}
\end{equation*}

If $\sigma_i|i|+\sigma_j|j|+\sigma_k|k|+\sigma_l|l|=0$, then using Lemma 4 of P\"{o}schel's article, we get
\begin{equation*}
|\delta|\geq \min_{m\in(0,\frac{1}{4})\cup(\frac{1}{4},\frac{41}{4})}\left\{\frac{m}{\sqrt{n(n+1)+m}},\frac{n}{\sqrt{m+2}},{\frac{4m-1}{4(n(n+1)+m)^{\frac{3}{2}}}} \right\}.
\end{equation*}

\textbf{Case 2.} Assume $\frac{1}{4}<m<\frac{41}{4}$, consider the following two cases
\begin{equation*}\left\{
\begin{aligned}
&i\pm j \pm k \pm l =2\alpha\quad \alpha\in \mathbb{Z}_+,\\
&i\pm j \pm k \pm l =2\beta\quad \beta\in \mathbb{Z}_-,
\end{aligned}\right.
\end{equation*}
where $0<i\leq j\leq k\leq l$.
The case ``$i\pm j\pm k\pm l =0$ '' is discussed in Lemma 4 of P\"{o}schel's article.

Before proving \textbf{Lemma 4.2} in \textbf{Case 2.}, we have to show another lemma and state a basic property of Legendre polynomials.
\begin{lem}
Assume $\frac{1}{4}<m<\frac{41}{4}$, $j-i>l-k,0< i\leq j\leq k\leq l$, and $i\pm j\pm k\pm l$ is even, if $\delta=(\lambda_j-\lambda_i)-(\lambda_l-\lambda_k)$, then $|\delta|\geq 2-\frac{4m-1}{4\sqrt{i(i+1)+m}+4i+2}>0$.
\end{lem}
\textbf{Proof.}\quad It is easy to get $j-i\geq l-k+2$, hence we have
\begin{equation*}
\begin{aligned}
|\delta|=&|(\lambda_j-\lambda_i)-(\lambda_l-\lambda_k)|\\
=&|[(\lambda_j-\lambda_i)-(l-k)]-[(\lambda_l-\lambda_k)-(l-k)]|\\
\geq &|(\lambda_j-\lambda_i)-(l-k)|-|(\lambda_l-\lambda_k)-(l-k)|\\
\geq & |\lambda_j-\lambda_i|-|l-k|+[(l-k)-(\lambda_l-\lambda_k)]\\
\geq & 2-[(j-i)-(\lambda_j-\lambda_i)]+[(l-k)-(\lambda_l-\lambda_k)]\\
\geq & 2-[(j-i)-(\lambda_j-\lambda_i)].
\end{aligned}
\end{equation*}
It is clear that
\begin{equation*}
\begin{aligned}
\lambda_j-\lambda_i=&\sqrt{j(j+1)+m}-\sqrt{i(i+1)+m}\\
=&\frac{(j+i+1)(j-i)}{\sqrt{j(j+1)+m}+\sqrt{i(i+1)+m}}\\
\leq&j-i,
\end{aligned}
\end{equation*}
when $\frac{1}{4}<m<\frac{41}{4}$.
Then
\begin{equation}
\begin{aligned}
&(j-i)-(\lambda_j-\lambda_i)\\
=&\left(1- \frac{(j+i+1)}{\sqrt{j(j+1)+m}+\sqrt{i(i+1)+m}} \right)(j-i)\\
=& \frac{\sqrt{j(j+1)+m}-(j+\frac{1}{2})}{\sqrt{j(j+1)+m}+\sqrt{i(i+1)+m}}\cdot(j-i)\\
+&\frac{\sqrt{i(i+1)+m}-(i+\frac{1}{2})}{\sqrt{j(j+1)+m}+\sqrt{i(i+1)+m}}\cdot(j-i)\\
=&\frac{(m-\frac{1}{4})(j-i)}{\left[\sqrt{j(j+1)+m}+(j+\frac{1}{2})\right]\left[\sqrt{j(j+1)+m}+\sqrt{i(i+1)+m}\right]}\\
+&\frac{(m-\frac{1}{4})(j-i)}{\left[\sqrt{i(i+1)+m}+(i+\frac{1}{2})\right]\left[\sqrt{j(j+1)+m}+\sqrt{i(i+1)+m}\right]}.\\
\end{aligned}
\end{equation}
If $j-i=h\geq 0$, then the function with respect to $h$ satisfies the following property,
\begin{equation*}
\begin{aligned}
f(h)=&(j-i)-(\lambda_j-\lambda_i)\\
=&h-\left[\sqrt{(i+h)(i+h+1)+m}-\sqrt{(i)(i+1)+m}\right]
\end{aligned}
\end{equation*}
is monotone increasing, since
\begin{equation*}
f^{\prime}(h)=1-\frac{2h+2i+1}{2\sqrt{(i+h)(i+h+1)+m}}>0,
\end{equation*}
when $\frac{1}{4}<m<\frac{41}{4}$.

Hence, we completes \textbf{Lemma 4.3.} by letting $h$ goes to infinity,
\begin{equation}\label{m}
\begin{aligned}
|\delta|\geq& 2-f(h)\\
\geq&2-\lim_{h\rightarrow\infty}f(h)\\
=&2-\frac{m-\frac{1}{4}}{\sqrt{i(i+1)+m}+i+\frac{1}{2}}>0.
\end{aligned}
\end{equation}
In fact, $\frac{1}{4}<m<\frac{41}{4}$ comes from the last inequality of (\ref{m}) if $i=1$  .
\qed\\

In order to prove \textbf{Case 2.}, we need to divide it into the following 9 subcases:\\

\textbf{Subcase 2.1.} $i+j+k+l=2\alpha$. Since $i+j+k+l\geq i+3j$, we have $\alpha>j$. Then, we get $i-j+k+l = 2(\alpha-j)$, which convert to Subcase 2.4 below.\\

\textbf{Subcase 2.2.} $i+j+k-l=2\alpha$. If $1\leq \alpha<i$, then we have $l-k=(i-\alpha)+(j-\alpha)>(j-\alpha)-(i-\alpha)=j-i$. Using the idea of \textbf{Lemma 4 in P\"{o}schel } \cite{poschel1996quasi} , one can obtain
\begin{equation*}
|\delta|\geq 2(i-\alpha)f^{\prime}(j)\geq \frac{2j+1}{\sqrt{j(j+1)+m}}.
\end{equation*}
On the other hand, it is easy to obtain $\alpha\leq \frac{i+j}{2}$. If $\alpha = i$ or $\alpha = \frac{i+j}{2}$, then it converts to P\"{o}schel's case. So it suffices to consider the case $i<\alpha<\frac{i+j}{2}$, which means $l-k<j-i$. This can be solved by using Lemma 4.3. \\

\textbf{Subcase 2.3.} $i+j-k+l=2\alpha$. Using the basic assumption, we can get $l-k\leq \alpha \leq i+j$. The case $\alpha = l-k$ or $\alpha =i+j$ can be solved by using \textbf{Lemma 4 in P\"{o}schel } \cite{poschel1996quasi} . If $\alpha>j$, then $i-j-k+l=2(\alpha-j)$, which converts to Subcase 2.5 below. So it suffices to consider the case $l-k<\alpha \leq j$, which means $l-k\leq j-i$. Use \textbf{Lemma 4 in P\"{o}schel } \cite{poschel1996quasi}  when the equality holds, while use Lemma 4.3 when equality does not hold.\\

\textbf{Subcase 2.4.} $i-j+k+l=2\alpha$. Using the basic assumption, we get $\alpha\geq l-j$. If $\alpha>k$, then $i-j-k+l=2(\alpha-k)$, which converts to Subcase 2.5 below. Otherwise $j-i=k+l-2\alpha\geq k+l -2k=l-k$, then use the same skill like Subcase 2.3. \\

\textbf{Subcase 2.5.} $i-j-k+l=2\alpha$. It is easy to be solved when we observe that $l-k=j-i+2\alpha$ by using the idea of  \textbf{Lemma 4 in P\"{o}schel } \cite{poschel1996quasi} .\\

\textbf{Subcase 2.6.} $i-j-k-l=2\beta$. If $\beta+l\geq 1$, then it converts to Subcase 2.5. Next, using the basic assumption, we get $-j-k\leq \beta \leq -l$, which means $i-j-k-l=2\beta\leq -2l$. This concludes that $l-k\leq j-i$, then we can use the same skill above.\\

\textbf{Subcase 2.7.} $i-j-k+l=2\beta$. Observe that $l-k=j-i-2|\beta|$.\\

\textbf{Subcase 2.8.} $i-j+k-l=2\beta$. Using the basic assumption, we get $|\beta|\leq k$. If $\beta+l<1$, then $i+j-k-l=2\beta\leq-2l$, i.e., $i+j+l-k\leq 0$, which is a contradiction. So we obtain $i+j-k+l=2(\beta+l)$, which can converts to Subcase 2.3.\\

\textbf{Subcase 2.9.} $i-j+k-l=2\beta$. Observe that $l-k+j-i=2|\beta|$, then $l-k\leq j-i+2$ or $l-k\geq j-i$, which can use the same skill as above.

Hence, we finish the proof of Lemma 4.2.
\qed

\section{The Cantor Manifold Theorem}
In this section, we will state Cantor manifold theorem in P\"{o}schel's  article \cite{poschel1996quasi} which is proven in \cite{kuksin1996invariant} using the KAM-theorem for partial differential equations from \cite{poschel2006kam}. The difficulty here is to check the nondegeneracy condition (\ref{non}) for Cantor manifold theorem.

In a neighbourhood of the origin in $\ell^2_s$, we now consider more generally hamiltonian of the form $H=\Lambda+Q+R$, where $\Lambda+Q$ is integrable and in normal form and $R$ is a perturbation term. Letting $z=(\tilde{z},\hat{z})$ with $\tilde{z}=(z_1,z_2)$ and $\hat{z}=(z_3,z_4,\ldots)$, as well as
\begin{equation*}
I=\frac{1}{2}(|z_1|^2,|z_2|^2),\quad Z=\frac{1}{2}(|z_3|^2,|z_4|^2,\ldots),
\end{equation*}
we assume that
\begin{equation*}
\Lambda=\langle \alpha,I\rangle+\langle \beta,Z\rangle,\quad Q=\langle AI,I\rangle+\langle BI,Z\rangle,
\end{equation*}
with constant vectors $\alpha,\beta$ and constant matrices $A,B$,
\begin{equation}\label{al1}
\alpha_i=\lambda_i (i=1,2),
\end{equation}
\begin{equation}\label{al2}
\beta_j=\lambda_j (j\geq 3).
\end{equation}
\begin{equation}\label{al3}
A=\left(
\begin{aligned}
&\bar{G}_{11}\quad\bar{G}_{12}\\
&\bar{G}_{21}\quad\bar{G}_{22}
\end{aligned}\right)
\end{equation}
\begin{equation}\label{al4}
B_{jk}=\bar{G}_{jk}(j\geq 3, k=1,2),
\end{equation}

In the Birkhoff normal form lemma, $\Lambda+\bar{G}$ is of that form.

The equations of motion of the hamiltonian $\Lambda+Q$ are
\begin{equation*}
\dot{\tilde{z}}=\mi(\alpha+AI+B^Tz)_j\tilde{z}_j,\quad \dot{\hat{z}}_j=\mi(\beta+BI)_j\hat{z}_j.
\end{equation*}
Thus, the complex $2$-dimensional manifold $E={\hat{z}=0}$ is invariant, and it is completely filled up to the origin by the invariant tori
\begin{equation*}
\mathscr{T}(I)=\{\tilde{z}:|\tilde{z}_j|^2=2I_j,1\leq j\leq 2\},\quad I\in\overline{\mathbb{P}^2}.
\end{equation*}
On $\mathscr{T}(I)$ the flow is given by the equations
\begin{equation}\label{omega}
\dot{\tilde{z}}_j=\mi\omega_j(I)\tilde{z}_j,\quad \omega(I)=\alpha+AI,
\end{equation}
and in its normal space by
\begin{equation}\label{Omega}
\dot{\hat{z}}_j=\mi\Omega_j(I)\hat{z}_j,\quad \Omega(I)=\beta+BI.
\end{equation}
They are linear and in diagonal form. In particular, since $\Omega(I)$ is real, $\hat{z}=0$ is an elliptic fixed point, all the tori are linearly stable, and their orbits have zero Lyapunov exponents. The Cantor manifold theorem proves the persistence of a large portion of $E$ forming an invariant Cantor manifold $\mathscr{E}$ for the hamiltonian $H=\Lambda+Q+R$.

For the existence of $\mathscr{E}$, the following assumptions are made.

\textbf{A.Nondegeneracy.} The normal form $\Lambda+Q$ is nondegenerate in the sense that
\begin{equation}\label{non}
\begin{aligned}
(A_1)\quad&\det{A}\neq 0,\\
(A_2)\quad&\langle l,\beta \rangle \neq 0,\\
(A_3)\quad&\langle k,\omega(I)\rangle+\langle l,\Omega(I)\rangle\not\equiv 0,
\end{aligned}
\end{equation}
for all $(k,l)\in \mathbb{Z}^2\times \mathbb{Z}^\infty$ with $1\leq|l|\leq2$.

\textbf{B.Spectral asymptotics.} There exists $d\geq 1$ and $\delta<1$ such that
\begin{equation*}
\beta_j=j^d+\ldots+O(j^\delta),
\end{equation*}
where the dots stand for terms of order less than $d$ in $j$. Note that the normalization of the coefficient of $j^d$ can always be achieved by a scaling of time.

\textbf{C.Regularity.}
\begin{equation*}
X_Q,X_R\in A(\ell^2_s,\ell^2_{\bar{s}}),\quad\left\{
\begin{aligned}
&\bar{s}\geq s \quad \text{for}\quad d>1,\\
&\bar{s}>s \quad \text{for}\quad d=1.
\end{aligned}\right.
\end{equation*}

By the regularity assumption, the coefficients of $B=(B_{ij})_{1\leq j \leq 2 < i}$ satisfy the estimate $B_{ij}=O(i^{s-\bar{s}})$ uniformly in $1\leq j \leq 2$. Consequently, for $d=1$ there exists a positive constant $\kappa$ such that
\begin{equation*}
\frac{\Omega_i-\Omega_j}{i-j}=1+O(j^{-\kappa}),\quad i>j,
\end{equation*}
uniformly for bounded $I$. For $d>1$, we set $\kappa=\infty$.

The following theorem is in P\"{o}schel \cite{poschel1996quasi}.

\begin{thm}
``\emph{THE CANTOR MANIFOLD THEOREM}. Suppose the hamiltonian $H=\Lambda+Q+R$ satisfies assumptions $A$,$B$ and $C$, and
\begin{equation*}
|R|=O(\norm{\hat{z}_s^4})+O(\norm{z}_s^g)
\end{equation*}
with
\begin{equation*}
g>4+\frac{4-\Delta}{\kappa},\quad \Delta=\min{(\bar{s}-s,1)}.
\end{equation*}
Then there exists a Cantor manifold $\mathscr{E}$ of real analytic, elliptic diophantine $n-$ tori given by a Lipschitz continuous embedding $\Psi:\mathscr{T}[\mathscr{C}]\rightarrow \mathscr{E}$, where $\mathscr{C}$ has full density at the origin, and $\Psi$ is close to the inclusion map $\Psi_0$:
\begin{equation*}
\norm{\Psi-\Psi_0}_{\bar{s},B_r\cap\mathscr{T}[\mathscr{C}]}=O(r^\sigma),
\end{equation*}
with some $\sigma>1$. Consequently, $\mathscr{E}$ is tangent to $E$ at the origin.''
\end{thm}

We now verify the assumptions of the Cantor Manifold Theorem. We already known that $X_Q,X_R\in A(\ell^2_{\frac{7}{2}},\ell^2_{\frac{9}{2}})$ with $|R|=O(\norm{\tilde{z}}^4_{\frac{7}{2}})+O(\norm{z}^6_{\frac{7}{2}})$. On the other hand, we have
\begin{equation*}
\lambda_j=\sqrt{j(j+1)+m}=j+\frac{1}{2}+\frac{m-\frac{1}{4}}{2j}+O(j^{-3}).
\end{equation*}
So conditions B and C are satisfied with $d=1,\delta=-1$,$\bar{s}=\frac{9}{2}$ and $s=\frac{7}{2}$.

Moreover, since $B_{ij}=\bar{G}_{ij}=\frac{2(2i+1)(2j+1)}{\lambda_i\lambda_j}P(i,j)$, we have
\begin{equation*}
\Omega_{j-2}=(\beta+BI)_{j-2}=\lambda_j+\frac{\langle v,I\rangle}{\lambda_j}
\end{equation*}
with $v=2(2i+1)(2j+1)P(i,j)(\lambda_1^{-1},\lambda_2^{-1})$. This gives the asymptotic expansion
\begin{equation*}
\Omega_{j-2}=j+\frac{1}{2}+\frac{m-\frac{1}{4}}{2j}+\frac{\langle v,I\rangle}{j}+O(j^{-3})=j+\frac{1}{2}+\frac{m_I}{j}+O(j^{-3}),
\end{equation*}
$m_I=\frac{m-\frac{1}{4}}{2}+\langle v,I\rangle$. Thus, for $i>j$,
\begin{equation*}
\frac{\Omega_i-\Omega_j}{i-j}=1-\frac{m_I}{(i+2)(j+2)}+O(j^{-3})=1+O(j^{-2}),
\end{equation*}
we can obtain
\begin{equation*}
\begin{aligned}
\det{A}&=\det{\left(
\begin{aligned}
&\bar{G}_{11}\quad\bar{G}_{12}\\
&\bar{G}_{21}\quad\bar{G}_{22}
\end{aligned}\right)}\\
&\triangleq\det{\left(
\begin{aligned}
&\frac{g_{11}}{\lambda_1^2}\quad&\frac{g_{12}}{\lambda_1\lambda_2}\\
&\frac{g_{21}}{\lambda_2\lambda_1}\quad&\frac{g_{22}}{\lambda_2^2}
\end{aligned}\right)}\\
&\triangleq\frac{1}{\lambda_1^2\lambda_2^2}\det{g}<0,
\end{aligned}\\
\end{equation*}
where
\begin{equation}
\begin{aligned}
&g_{11}=\lambda_1^2\bar{G}_{11}=3\times 3 \times \mathbf{P}(1,1)=\frac{18}{5},\\
&g_{12}=g_{21}=\lambda_1\lambda_2\bar{G}_{12}=2\times 3 \times 7 \times \mathbf{P}(1,3)=\frac{92}{15}, \\
&g_{22}=\lambda_2^2\bar{G}_{22}=7\times 7 \times \mathbf{P}(3,3)=\frac{2\times 7\times 241}{5\times 11\times 13},
\end{aligned}
\end{equation}
\begin{equation}\label{deg}
\deg{g}=g_{11}g_{22}-g_{12}g_{21}=-\frac{663764}{32175}.
\end{equation}
The nondegeneracy condition $(A_2)$ is easy to check since $\lambda_j$ or $\lambda_i\pm \lambda_j(i\neq j)$ are not equal to zero.

Next, we will check the nondegeneracy condition $(A_3)$. Since the condition $1\leq|l|\leq2$, we only need to consider the following two cases£º
\begin{equation}\label{kw1}
\langle k,\omega \rangle\pm\Omega_j\not\equiv 0
\end{equation}
and
\begin{equation}\label{kw2}
\langle k,\omega \rangle\pm(\Omega_i-\Omega_j)\not\equiv 0
\end{equation}
Recall the definition of $\omega(I)$ (\ref{omega}) and  $\Omega(I)$ (\ref{Omega}), we can obtain that $\Omega=\beta+BA^{-1}(\omega-\alpha)$. Besides, choosing $\zeta=(\frac{1}{\omega_1},\frac{\omega_2}{\omega_1})=(\sigma,\zeta_2)$ as another new parameter vector instead of $\omega=(\omega_1,\omega_2)$, we can obtain the following two expressions $(A_{31})$ and $(A_{32})$  with respect to $\sigma=\frac{1}{\omega_1}$ which are equivalent to (\ref{kw1}) and (\ref{kw2}).
\begin{equation*}
(A_{31})\quad \left|\frac{\md f_1(\sigma)}{\md \sigma} \right|=\left|\beta_j-\sum_{l=1}^{2}B_{jl}(A^{-1}\alpha)_l\right|>0,
\end{equation*}

\begin{equation*}
(A_{32})\quad \left|\frac{\md f_2(\sigma)}{\md \sigma} \right|=\left|\beta_i-\beta_j-\sum_{l=1}^{2}(B_{il}-B_{jl})(A^{-1}\alpha)_l\right|>0,
\end{equation*}
 where
 \begin{equation}
 \tilde{\omega}=\omega\sigma=(1,\zeta_2),
 \end{equation}
 \begin{equation}
 f_1(\sigma)=\langle k,\tilde{\omega} \rangle+\sigma\beta_j+\sum_{l=1}^{2}B_{jl}[A^{-1}(\tilde{\omega}-\alpha\sigma)]_l,
 \end{equation}
 \begin{equation}
 f_2(\sigma)=\langle k,\tilde{\omega} \rangle\pm\left[\left(\sigma\beta_i+\sum_{l=1}^{2}B_{il}[A^{-1}(\tilde{\omega}-\alpha\sigma)]_l\right)-\left(\sigma\beta_j+\sum_{l=1}^{2}B_{jl}[A^{-1}(\tilde{\omega}-\alpha\sigma)]_l\right)\right].
\end{equation}

Here, we denote $g_{j1}$ and $g_{j2}$ according to (\ref{G2}),
\begin{equation}\label{dg1}
\begin{aligned}
g_{j1}&=\lambda_1\lambda_j\bar{G}_{1j}\\
&=2\times (4\times 1-1)\times (4j-1)\times P(1,2j-1)\\
&=\frac{12(8j^2-4j-1)}{(4j-3)(4j+1)},
\end{aligned}
\end{equation}
and
\begin{equation}\label{dg2}
\begin{aligned}
g_{j2}&=\lambda_2\lambda_j\bar{G}_{2j}\\
&=2\times(4\times 2-1)\times(4j-1)P(3,2j-1)\\
&=\frac{28\times(1088j^6-1632j^5-2440j^4+3120j^3+1406j^2-1110j-225)}{(4j-7)(4j-5)(4j-3)(4j+1)(4j+3)(4j+5)}.
\end{aligned}
\end{equation}
By the basic computation
\begin{equation}
\begin{aligned}
A^{-1}&=\left(
\begin{aligned}
&\bar{G}_{11}\quad\bar{G}_{12}\\
&\bar{G}_{21}\quad\bar{G}_{22}
\end{aligned}\right)^{-1}\\
&=\frac{\left(
\begin{aligned}
&\bar{G}_{22}\quad-\bar{G}_{12}\\
&-\bar{G}_{21}\quad\bar{G}_{11}
\end{aligned}
\right)
}{\det{A}},
\end{aligned}
\end{equation}
as well as the definition of $B_{jk}$, $\alpha_i$, $\beta_j$ and $g_{j1}, g_{j2}$ in (\ref{al1})---(\ref{al4}) and (\ref{dg1}) (\ref{dg2}).
we can obtain the following result.\\

In $(A_{31})$, we have
\begin{equation}\label{f1}
\begin{aligned}
\left|\frac{\md f_1(\sigma)}{\md \sigma}\right|
&=\left| \lambda_j-\bar{G}_{j1}\frac{\bar{G}_{22}\lambda_1-\bar{G}_{12}\lambda_2}{\det{A}}- \bar{G}_{j2}\frac{\bar{G}_{11}\lambda_2-\bar{G}_{12}\lambda_1}{\det{A}} \right|\\
&=\frac{1}{|\det{A}|}\left|(\det{A})\lambda_j-\frac{g_{j1}}{\lambda_j\lambda_1}\left(\frac{g_{22}}{\lambda_2^2}\lambda_1-\frac{g_{12}}{\lambda_1\lambda_2}\lambda_2 \right)-\frac{g_{j2}}{\lambda_j\lambda_2}\left(\frac{g_{11}}{\lambda_1^2}\lambda_2-\frac{g_{12}}{\lambda_1\lambda_2}\lambda_1 \right) \right|\\
&=\frac{1}{|\det{A}|}\frac{1}{\lambda_1^2\lambda_2^2\lambda_j}\left|(\det{g})\lambda_j^2-g_{j1}(g_{22}\lambda_1^2-g_{12}\lambda_2^2)-g_{j2}(g_{11}\lambda_2^2-g_{12}\lambda_1^2) \right|\\
&\geq\frac{1}{|\det{A}|}\frac{1}{\lambda_1^2\lambda_2^2\lambda_j}\left[\left|(\det{g})\lambda_j^2\right|-\left|g_{j1}(g_{22}\lambda_1^2-g_{12}\lambda_2^2)+g_{j2}(g_{11}\lambda_2^2-g_{12}\lambda_1^2) \right|\right],
\end{aligned}
\end{equation}
where $j\geq 3$ and $m\in(0,\frac{1}{4})\cup(\frac{1}{4},\frac{41}{4})$. An easy computation gives
\begin{equation}\label{gg1}
  g_{22}\lambda_1^2-g_{12}\lambda_2^2=-\frac{3034}{2145}m-\frac{45876}{715}.
\end{equation}
and
\begin{equation}\label{gg2}
g_{11}\lambda_2^2-g_{12}\lambda_1^2=-\frac{38}{15}m+\frac{464}{15}.
\end{equation}

From (\ref{deg}) (\ref{gg1}) (\ref{gg2}), it is easy to estimate that
\begin{equation}\label{f11}
|(\det{g})\lambda_j^2|=|(\det{g})\left(2j(2j-1)+m\right)|>618,
\end{equation}
and
\begin{equation}\label{f12}
\left|g_{j1}(g_{22}\lambda_1^2-g_{12}\lambda_2^2)+g_{j2}(g_{11}\lambda_2^2-g_{12}\lambda_1^2) \right|<441.
\end{equation}
Then, by substituting (\ref{f11}) and (\ref{f12}) into (\ref{f1}), we obtain
\begin{equation*}
\left|\frac{\md f_1(\sigma)}{\md \sigma}\right|>0.
\end{equation*}
In $(A_{32})$, we have
\begin{equation}\label{f2}
\begin{aligned}
\left|\frac{\md f_2(\sigma)}{\md \sigma}\right|&=\left|\beta_i-\beta_j-\sum_l(B_{il}-B_{jl})(A^{-1}\alpha)_l \right|\\
=&\left| (\lambda_i-\lambda_j)-(\bar{G}_{i1}-\bar{G}_{j1})\frac{\bar{G}_{22}\lambda_1-\bar{G}_{12}\lambda_2}{\det{A}}
-(\bar{G}_{i2}-\bar{G}_{j2})\frac{\bar{G}_{11}\lambda_2-\bar{G}_{12}\lambda_1}{\det{A}} \right|\\
=&\frac{1}{\lambda_1^2\lambda_2^2|\det{A}|}\left| (\lambda_i-\lambda_j)\det{g}-\lambda_1^2\lambda_2^2\left(\frac{g_{i1}}{\lambda_i\lambda_1}-\frac{g_{j1}}{\lambda_j\lambda_1}\right)\left(\frac{g_{22}}{\lambda_2^2}\lambda_1-\frac{g_{12}}{\lambda_1\lambda_2}\lambda_2\right)\right.\\
&\left.-\lambda_1^2\lambda_2^2\left(\frac{g_{i2}}{\lambda_i\lambda_2}-\frac{g_{j2}}{\lambda_j\lambda_2}\right)\left(\frac{g_{11}}{\lambda_1^2}\lambda_2-\frac{g_{12}}{\lambda_1\lambda_2}\lambda_1\right) \right|\\
\geq&\frac{|i-j|}{\lambda_1^2\lambda_2^2|\det{A}|(\lambda_i+\lambda_j)}\left[\left| (\det{g})\frac{\lambda_i-\lambda_j}{i-j}\right| -\left|(g_{22}\lambda_1^2-g_{12}\lambda_2^2)\frac{\frac{g_{i1}}{\lambda_i}-\frac{g_{j1}}{\lambda_j}}{i-j}\right.\right.\\
&+\left.\left.(g_{11}\lambda_2^2-g_{12}\lambda_1^2)\frac{\frac{g_{i2}}{\lambda_i}-\frac{g_{j2}}{\lambda_j}}{i-j}\right|\right].
\end{aligned}
\end{equation}
Since
\begin{equation*}
\begin{aligned}
\frac{\lambda_i-\lambda_j}{i-j}&=\frac{\sqrt{2i(2i-1)+m}-\sqrt{2j(2j-1)+m}}{i-j}\\
&=\frac{4(i+j)-2}{\lambda_i+\lambda_j}\\
&>\frac{4(i+j)-2}{[2(i+j)-1]+2\sqrt{m-\frac{1}{4}}}\\
&>\frac{2}{1+\frac{\sqrt{40}}{13}},\quad(i,j\geq 3,i\not=j)
\end{aligned}
\end{equation*}
we have the following estimate by (\ref{deg})
\begin{equation}\label{f21}
\left| (\det{g})\frac{\lambda_i-\lambda_j}{i-j}\right|>27.
\end{equation}
Let
\begin{equation}\label{v1}
  v_1(i)=\frac{g_{i1}}{\lambda_i},
\end{equation}
\begin{equation}\label{v2}
  v_2(i)=\frac{g_{i2}}{\lambda_i}.
\end{equation}
Then for every fixed $m\in(0,\frac{1}{4})\cup (\frac{1}{4},\frac{41}{4})$, when $x\in[3,\infty)$
\begin{equation}
  \partial_xv_1(x)<0, \quad \partial_xv_2<0,
\end{equation}
and both monotonically increase with respect to $x$. Besides,
\begin{equation}
  \frac{\frac{g_{i1}}{\lambda_i}-\frac{g_{j1}}{\lambda_j}}{i-j}\in\partial_xv_1([3,\infty)),\quad \frac{\frac{g_{i2}}{\lambda_i}-\frac{g_{j1}}{\lambda_j}}{i-j}\in\partial_xv_2([3,\infty)).
\end{equation}
From (\ref{gg1}) (\ref{gg2}) and some computation, we know
\begin{equation}\label{f22}
\max_{m\in(0,\frac{1}{4})\cup (\frac{1}{4},\frac{41}{4})}\left|(g_{22}\lambda_1^2-g_{12}\lambda_2^2)\frac{\frac{g_{i1}}{\lambda_i}-\frac{g_{j1}}{\lambda_j}}{i-j}+(g_{11}\lambda_2^2-g_{12}\lambda_1^2)\frac{\frac{g_{i2}}{\lambda_i}-\frac{g_{j2}}{\lambda_j}}{i-j}\right|<27.
\end{equation}
In the end, by substituting (\ref{f21}) and (\ref{f22}) into (\ref{f2}), we  get
\begin{equation*}
\left|\frac{\md f_2(\sigma)}{\md \sigma}\right|>0.
\end{equation*}
\qed

Thus the main theorem follows.

\section{Proof of Theorem 3.2.}
Before our proof, we need the following four lemmas which will be proved in Appendix B.

 The idea is using the recursion formula (\ref{recursion}) to obtain lemma \ref{lp2} (with the help of lemma \ref{lp1} also ) and lemma \ref{lp3}. Then, we eliminate $\mathbf{Q}(m,n)$ and get the recursion (See lemma \ref{lp4}.) with respect to $\mathbf{P}(m,n)$ by some technical calculation.

The estimate of $\mathbf{c}(m)(=\lim_{n\rightarrow \infty}n\mathbf{P}(m,n))$ in lemma \ref{lp5} and the symmetry of $\mathbf{P}(m,n)$ give the final result. The idea using another sequence $\frac{O(1)}{m+1}$ similar to $\mathbf{c}(m)$ to obtain a rough upper bound of $\mathbf{c}(m)$ and using the symmetry property to obtain a precise estimate.\\

\begin{lem}\label{lp1}
$\int_{-1}^{1}P_m^2P_{n+1}P_{n-1}\md x=\int_{-1}^{1}P_n^2P_{m+1}P_{m-1}\md x$.
\end{lem}

If we denote
\begin{equation*}
\begin{aligned}
&\mathbf{P}(m,n)=\int_{-1}^{1}P_mP_mP_nP_n\md x,\\
&\mathbf{Q}(m,n)=\int_{-1}^{1}P_{m-1}P_{m+1}P_nP_n\md x,
\end{aligned}
\end{equation*}
then we have
\begin{lem}\label{lp2}
\begin{equation}\label{P}
  \begin{aligned}
  \mathbf{P}(m,n)=&\left(\frac{n+1}{2n+1}\right)^2\left(\frac{2m-1}{m}\right)^2\mathbf{P}(m-1,n-1)\\
  &+\frac{2n^2+2n-2m^2+2m}{(2n+1)^2m^2}\mathbf{Q}(m-1,n)-\left(\frac{m-1}{m}\right)^2\mathbf{P}(m-2,n).
  \end{aligned}
\end{equation}
\end{lem}

\begin{lem}\label{lp3}
\begin{equation}\label{Q}
\begin{aligned}
  \mathbf{Q}(m,n)=&\frac{(2m+1)m}{(m+1)(2m-1)}\mathbf{P}(m,n)\\
  &+\frac{(2m+1)(m-1)}{(2m-1)(m+1)}\mathbf{Q}(m-1,n)-\frac{m}{m+1}\mathbf{P}(m-1,n).
\end{aligned}
\end{equation}
\end{lem}

\begin{lem}\label{lp4}
\begin{equation}\label{PP}
\begin{aligned}
\mathbf{P}(m+1,n)=&\alpha_m^{n-1}\mathbf{P}(m,n-1)+\alpha_m^n\mathbf{P}(m,n)+\alpha_m^{n+1}\mathbf{P}(m,n+1)\\
-&\alpha_{m-1}^{n-1}\mathbf{P}(m-1,n-1)-\alpha_{m-1}^n\mathbf{P}(m-1,n)-\alpha_{m-1}^{n+1}\mathbf{P}(m-1,n+1) \\
+&\alpha_{m-2}^n\mathbf{P}(m-2,n),
\end{aligned}
\end{equation}
with the coefficients $\alpha_i^j$ and $\mathbf{P}(k,n),(k=0,1,2,3)$ expressed in Theorem 3.3.
\end{lem}

Now we only need to prove
\begin{equation*}
\int_{-1}^{1}P_i(x)P_j(x)P_k(x)P_l(x)\md x\lesssim \frac{1}{\sqrt{i+\frac{1}{2}}\sqrt{j+\frac{1}{2}}\sqrt{k+\frac{1}{2}}\sqrt{l+\frac{1}{2}}}.
 \end{equation*}

It is obvious to conclude the following proposition by mathematical induction according to the recursion (\ref{PP}) and the expressions of $\textbf{P}(i,n),(i=0,1,2)$ (\ref{PN0}) (\ref{PN1}) (\ref{PN2}).
\begin{prop}
For every fixed $m$, $\mathbf{P}(m,n)$ is a rational fraction with respect to $n$, i.e.
\begin{equation*}
  \mathbf{P}(m,n)=\frac{\mathscr{R}(n)}{\mathscr{S}(n)}.
\end{equation*}
The degree of denominator $\deg(\mathscr{R}(n))$ and the degree of numerator $\deg(\mathscr{S}(n))$ satisfy
\begin{equation*}
  \deg(\mathscr{S}(n))-\deg(\mathscr{R}(n))\geq1.
\end{equation*}
\end{prop}
\begin{rem}
  We conjecture that $\deg(\mathscr{R}(n))=2m$ and  $\deg(\mathscr{S}(n))=2m+1$.
\end{rem}
An immediate consequence of the proposition is that, for every fixed $m$,
\begin{equation*}
  \lim_{n\rightarrow \infty} n\mathbf{P}(m,n)= O(1),
\end{equation*}
and since the symmetry $\mathbf{P}(m,n)=\mathbf{P}(n,m)$, we also have that for every fixed $n$,
\begin{equation*}
  \lim_{m\rightarrow \infty}m\mathbf{P}(m,n)= O(1).
\end{equation*}

In the following, we have to check if,
\begin{equation}
  \lim_{mn\rightarrow \infty}mn\mathbf{P}(m,n)= O(1),
\end{equation}
which seems to be correct intuitively but need to be proved strictly.

If we fix $m$, we can denote
\begin{equation}\label{PMN}
  \lim_{n\rightarrow \infty}n\mathbf{P}(m,n)=\mathbf{c}(m),
\end{equation}
and it is obvious that
\begin{equation*}
  \mathbf{c}(0)=1,\quad\mathbf{c}(1)=\frac{1}{2},\quad\mathbf{c}(2)=\frac{11}{32},\quad\mathbf{c}(3)=\frac{17}{64}.
\end{equation*}
Observing the relationship
\begin{equation*}
  \alpha_m^{n-1}+\alpha_m^n+\alpha_m^{n+1}-\alpha_{m-1}^{n-1}-\alpha_{m-1}^n-\alpha_{m-1}^{n+1}+\alpha_{m-2}^n \equiv 1,
\end{equation*}
by letting $n$ goes to infinity, we get
\begin{equation}\label{idd}
  \alpha_m-\alpha_{m-1}+\alpha_{m-2}\equiv 1,
\end{equation}
where
\begin{equation}
  \alpha_i=\lim_{n\rightarrow \infty}\sum_{j=n-1}^{n+1}\alpha_i^j,(i=m-1,m),\quad \alpha_{m-2}=\lim_{n\rightarrow \infty}\alpha_{m-2}^n,
\end{equation}
and
\begin{equation}
  \alpha_{m-1}=\frac{6m^3-2m^2+1}{2(m+1)^3},\quad \alpha_m=\frac{(2m+1)(6m^3+2m^2-1)}{2(2m-1)(m+1)^3},\quad \alpha_{m-2}=\frac{(m-1)^3(2m+1)}{(m+1)^3(2m-1)}.
\end{equation}
Then from (\ref{PP}), letting $n$ goes to infinity, we have the following \emph{unilateral} sequence
\begin{equation}
\mathbf{c}(m+1)=\alpha_m \mathbf{c}(m)-\alpha_{m-1}\mathbf{c}(m-1)+\alpha_{m-2}\mathbf{c}(m-2).
\end{equation}
\begin{lem}\label{lp5}
The sequence $\mathbf{c}(m)$ satisfies the following properties:
\begin{equation}
0\leq \mathbf{c}(m+1)<\mathbf{c}(m)\leq \frac{1}{2}.
\end{equation}
and
\begin{equation}\label{upper}
\mathbf{c}(m+1)\leq\frac{15}{16(m+2)}+\frac{1}{32}.
\end{equation}
\end{lem}
The proof is given in \textbf{Appendix B}.

From (\ref{PMN}) and (\ref{upper}), we know
\begin{equation}\label{P1}
\sup_{n\geq 1} \left(n\mathbf{P}(m,n)\right)\lesssim \frac{15}{16(m+2)}+\frac{1}{32}.
\end{equation}

On the other hand, since the symmetry (See the definition of $\mathbf{P}(m,n)$. )
\begin{equation}
\mathbf{P}(m,n)=\mathbf{P}(n,m),
\end{equation}
\begin{equation}
\lim_{m\rightarrow \infty}m\mathbf{P}(m,n)= \lim_{m\rightarrow \infty}m\mathbf{P}(n,m)=\mathbf{c}(n),
\end{equation}
then by (\ref{upper})
\begin{equation}
m\mathbf{P}(m,n)\lesssim \mathbf{c}(n)=\frac{15}{16(n+1)}+\frac{1}{32},
\end{equation}
so
\begin{equation}\label{P2}
n\mathbf{P}(m,n)\lesssim \left(\frac{15n}{16(n+1)}+\frac{1}{32}n\right)\frac{1}{m}
\end{equation}

However,~from (\ref{P1}) we know
\begin{equation*}
\sup_{n\geq 1} \left(n\mathbf{P}(m,n)\right)\quad <\infty,
\end{equation*}
so from (\ref{P2}), we obtain
\begin{equation*}
\sup_{n\geq 1} \left(n\mathbf{P}(m,n)\right)\lesssim \frac{1}{m},
\end{equation*}
then it follows that, there exists $\widetilde{C}>0$ such that
\begin{equation*}
  \mathbf{P}(m,n)<\frac{\widetilde{C}}{mn}.
\end{equation*}
By Cauchy Inequality, we have
\begin{equation*}
\begin{aligned}
  \int_{-1}^{1}P_iP_jP_kP_l\md x\leq &\left(\int_{-1}^{1}(P_iP_j)^2\md x\right)^{\frac{1}{2}}\left(\int_{-1}^{1}(P_kP_l)^2\md x\right)^{\frac{1}{2}}\\
  \lesssim &\frac{1}{\sqrt{i+\frac{1}{2}}\sqrt{j+\frac{1}{2}}\sqrt{k+\frac{1}{2}}\sqrt{l+\frac{1}{2}}}.
\end{aligned}
\end{equation*}
If $i\leq j\leq k \leq l$, then we already know
\begin{equation*}
  \int_{-1}^{1}P_iP_jP_kP_l\md x\not=0 \quad(\text{iff}\quad i+j+k+l=2\alpha \quad \text{and} \quad i+j+k\geq l)
\end{equation*}
and $P_iP_jP_kP_l$ is an even function, from (\ref{Adam}) we know
\begin{equation*}
  \int_{-1}^{1}P_iP_jP_kP_l\md x>0,
\end{equation*}
then
\begin{equation*}
  0<\int_{-1}^{1}\phi_i\phi_j\phi_k\phi_l \md x\leq C
\end{equation*}
is an easy deduction.
\qed

\section{Appendix A}

Recall the definition
\begin{equation}\label{mixnorm}
  \{f\}_\alpha=\left\{
  \begin{aligned}
  &[f]_\alpha=\sup_{x\neq y}\frac{|f(x)-f(y)|}{|x-y|^\alpha}\quad when \quad 0<\alpha<1,\\
  &\norm{f}_{L^p[-1,1]}\quad when \quad \alpha=-\frac{1}{p}\leq 0,
  \end{aligned}\right.
\end{equation}
for any $f$ defined on $[-1,1]$.

Then we have the following General Nirenberg Inequality.
\begin{thm}\cite{henry1982remember}\cite{nirenberg1966extended} \cite{nirenberg2011elliptic}
Suppose $-1\leq \alpha,\beta,\gamma <1$ and $j,k$ are nonnegative integers, $0\leq \theta \leq 1$, and
\begin{equation}
\left\{
  \begin{aligned}
  j+\beta&=\theta(k+\alpha)+(1-\theta)\gamma,\\
  \frac{1}{p_\beta}&\leq\frac{\theta}{p_\alpha}+\frac{1-\theta}{p_\gamma},
  \end{aligned}
\right.
\end{equation}
where
\begin{equation}
p_\delta=\left\{
  \begin{aligned}
  &\infty,\quad \delta\geq 0,\\
  &-\frac{1}{\delta},\quad \delta<0.
  \end{aligned}
\right.
\end{equation}
When $k+\alpha$ is an integer $\geq 1$, and $-1<\alpha<0$ (i.e. $1<p_\alpha<\infty$) we require $\theta\neq 1$. Then for any $f\in C^\infty_0[-1,1]$, we have
\begin{equation}\label{gn1}
  \{\partial_x^j f\}_\beta\lesssim\{\partial_x^k f\}_\alpha^\theta\{ f\}_\gamma^{1-\theta}.
\end{equation}
\end{thm}
An immediate consequence of this theorem is the Gagliardo-Nirenberg Inequality.
\begin{thm}\cite{nirenberg1966extended}\cite{nirenberg2011elliptic}
Suppose $0\leq a \leq 1,1\leq p,q,r\leq \infty,m,k$ are any integers satisfying
\begin{equation}
\left\{
  \begin{aligned}
  &k-\frac{1}{p}=a(m-\frac{1}{q})+(1-a)(-\frac{1}{r}),\\
  &\frac{1}{p}\leq\frac{a}{q}+\frac{1-a}{r},
  \end{aligned}
\right.
\end{equation}
then for any $f\in C^\infty_0[-1,1],$ we have the following inequality
\begin{equation}\label{gnin}
  \norm{\partial_x^ku}_{L^p[-1,1]}\lesssim\norm{\partial_x^mu}_{L^q[-1,1]}^a\norm{u}_{L^r[-1,1]}^{1-a}
\end{equation}
with the following exception: if $m-\frac{n}{q}=k,1<p<\infty$, then (\ref{gnin}) holds for $a\neq 1$.
\end{thm}
\section{Appendix B}

\subsection{Proof of Lemma {\ref{lp1}}.} From \cite{hs1938certain}, we know that
\begin{equation*}
\begin{aligned}
  &P_{n+1}P_{n-1}=\sum_{k=2}^{2n}\left(\begin{array}{clcr} n+1 & n-1 & k\\ 0 & \quad 0 & 0\end{array}\right)^2(2k+1)P_k\\
  &P_m^2=\sum_{k=0}^{2m}\binom{m \quad m\quad  k}{0\quad 0\quad  0}^2(2k+1)P_k,
\end{aligned}
\end{equation*}
and the same to $P_{m+1}P_{m-1}$ and $P_n^2$, then
\begin{equation*}
\begin{aligned}
  &\int_{-1}^{1}P_{n+1}P_{n-1}P_{m}^2\md x=2\sum_{k=2}^{2\min(m,n)}\left(\begin{array}{clcr} n+1 & n-1 & k\\ 0 & \quad 0 & 0\end{array}\right)^2\binom{m \quad m\quad  k}{0\quad 0\quad  0}^2(2k+1),\\
  &\int_{-1}^{1}P_{m+1}P_{m-1}P_{n}^2\md x=2\sum_{k=2}^{2\min(m,n)}\left(\begin{array}{clcr} m+1 & m-1 & k\\ 0 & \quad 0 & 0\end{array}\right)^2\binom{n \quad n\quad  k}{0\quad 0\quad  0}^2(2k+1).
\end{aligned}
\end{equation*}
It is easy to check that
\begin{equation*}
\begin{aligned}
&\left(\begin{array}{clcr} n+1 & n-1 & k\\ 0 & \quad 0 & 0\end{array}\right)^2\binom{m \quad m\quad  k}{0\quad 0\quad  0}^2\\
=&\frac{\mathbf{A}(\frac{k}{2}-1)\mathbf{A}(\frac{k}{2}+1)\mathbf{A}(n-\frac{k}{2})}{(2n+k+1)\mathbf{A}(n+\frac{k}{2})}
\frac{\mathbf{A}(\frac{k}{2})\mathbf{A}(\frac{k}{2})\mathbf{A}(m-\frac{k}{2})}{(2m+k+1)\mathbf{A}(m+\frac{k}{2})}\\
=&\frac{\mathbf{A}(\frac{k}{2}-1)\mathbf{A}(\frac{k}{2}+1)\mathbf{A}(m-\frac{k}{2})}{(2m+k+1)\mathbf{A}(m+\frac{k}{2})}
\frac{\mathbf{A}(\frac{k}{2})\mathbf{A}(\frac{k}{2})\mathbf{A}(n-\frac{k}{2})}{(2n+k+1)\mathbf{A}(n+\frac{k}{2})}\\
=&\left(\begin{array}{clcr} m+1 & m-1 & k\\ 0 & \quad 0 & 0\end{array}\right)^2\binom{n \quad n\quad  k}{0\quad 0\quad  0}^2.
\end{aligned}
\end{equation*}
\qed\\
\subsection{Proof of Lemma \ref{lp2}.}
From the recursion formula (\ref{recursion})
\begin{equation*}
(n+1)P_{n+1}-(2n+1)xP_n+nP_{n-1}=0,
\end{equation*}
we get
\begin{equation*}
  P_m=\frac{(2m-1)xP_{m-1}-(m-1)P_{m-2}}{m},\quad xP_n=\frac{(n+1)P_{n+1}+nP_{n-1}}{2n+1},
\end{equation*}
then we can write the demanded integration as follows,
\begin{equation}\label{e1}
  \begin{aligned}
  &\int_{-1}^{1}P_mP_mP_nP_n\md x=\int_{-1}^{1}\left[\frac{(2m-1)xP_{m-1}-(m-1)P_{m-2}}{m}\right]^2P_n^2\md x\\
  =&\left(\frac{2m-1}{m}\right)^2\int_{-1}^{1}P_{m-1}^2\left[\frac{n+1}{2n+1}P_{n+1}+\frac{n}{2n+1}P_{n-1}\right]^2\md x+\left(\frac{m-1}{m}\right)^2\int_{-1}^{1}P_{m-2}^2P_n^2\md x\\
  &-\frac{2(2m-1)(m-1)}{m^2}\left[\frac{mP_m+(m-1)P_{m-2}}{2m-1}\right]P_{m-2}P_n^2\md x\\
  =&\left(\frac{n+1}{2n+1}\right)^2\left(\frac{2m-1}{m}\right)^2\int_{-1}^{1}P_{m-1}^2P_{n+1}^2\md x+\left(\frac{2m-1}{m}\right)^2\left(\frac{n}{2n+1}\right)^2\int_{-1}^{1}P_{m-1}^2P_{n-1}^2\md x\\
  &+2\left(\frac{2m-1}{m}\right)^2\frac{n+1}{2n+1}\frac{n}{2n+1}\int_{-1}^{1}P_{m-1}^2P_{n+1}P_{n-1}\md x-\left(\frac{m-1}{m}\right)^2\int_{-1}^{1}P_{m-2}^2P_n^2\md x\\
  &-\frac{2(m-1)}{m}\int_{-1}^{1}P_mP_{m-2}P_n^2\md x.
  \end{aligned}
\end{equation}

With the help of\textbf{ Lemma \ref{lp1}}, we have
\begin{equation}
  \int_{-1}^{1}P_{m-1}^2P_{n+1}P_{n-1}\md x = \int_{-1}^{1}P_{m}P_{m-2}P_{n}^2\md x.
\end{equation}
Insert it into equation (\ref{e1}), we complete the proof.\\
\qed\\
\subsection{Proof of Lemma \ref{lp3}. }From (\ref{recursion}) we know
\begin{equation*}
  \begin{aligned}
  &\int_{-1}^{1}[(n+1)P_{n+1}]P_{n+1}P_m^2\md x\\
  =&\int_{-1}^{1}\left[(2n+1)xP_n\right]P_{n+1}P_m^2\md x-\int_{-1}^{1}\left[nP_{n-1}\right]P_{n+1}P_m^2\md x\\
  =&(2n+1)\int_{-1}^{1}(xP_{n+1})P_nP_m^2\md x-n\int_{-1}^{1}P_{n+1}P_{n-1}P_m^2\md x\\
  =&(2n+1)\int_{-1}^{1}P_nP_m^2\frac{(n+2)P_{n+2}+(n+1)P_n}{2n+3}\md x-n\int_{-1}^{1}P_{n-1}P_{n+1}P_m^2,
  \end{aligned}
\end{equation*}
\qed
\subsection{Proof of lemma \ref{lp4}.}
Using the recursion (\ref{recursion})
\begin{equation*}
  (n+1)P_{n+1}-(2n+1)xP_n+nP_{n-1}=0,
\end{equation*}
and the expressions of Legendre polynomials
\begin{equation*}
  P_0(x)=1,\quad P_1(x)=x,\quad P_2(x)=\frac{3x^2-1}{2},
\end{equation*}
it is easy to prove
\begin{equation*}
  \mathbf{P}(1,n)=\frac{2(2n^2+2n-1)}{(2n-1)(2n+1)(2n+3)},
\end{equation*}
and
\begin{equation*}
   \mathbf{P}(0,n)=\frac{2}{2n+1},
\end{equation*}
as well as
\begin{equation*}
\mathbf{P}(2,n)=\frac{11n^4+22n^3-31n^2-42n+18}{(2n-3)(2n-1)(2n+1)(2n+3)(2n+5)}.
\end{equation*}
Then
\begin{equation*}
  \mathbf{Q}(1,n)=\int_{-1}^{1}P_0P_1P_n^2\md x =\frac{3}{2}\mathbf{P}(1,n)-\frac{1}{2}\mathbf{P}(0,n)=\frac{2n(n+1)}{(2n-1)(2n+1)(2n+3)}.
\end{equation*}

From (\ref{P}) and (\ref{Q}), we can find the expressions of $\mathbf{P}(m,n)$ and $\mathbf{Q}(m,n)$, when $m$ is fixed. In particular, we have
\begin{equation*}
  \mathbf{P}(3,n)=\frac{34n^6+102n^5-305n^4-780n^3+703n^2+1110n-450}{(2n-5)(2n-3)(2n-1)(2n+1)(2n+3)(2n+5)(2n+7)}.
\end{equation*}\\

Now, by\textbf{ Lemma \ref{lp3}}, we have the following recursion with respect to $\{\frac{m(m+1)}{2m+1}\mathbf{Q}(m,n)\}_{m=1}^{\infty}$,
\begin{equation*}
  \frac{m(m+1)}{2m+1}\mathbf{Q}(m,n)-\frac{(m-1)[(m-1)+1]}{2(m-1)+1}\mathbf{Q}(m-1,n)=\frac{m^2}{2m-1}\mathbf{P}(m,n)-\frac{m^2}{2m+1}\mathbf{P}(m-1,n).
\end{equation*}
Then it follows that
\begin{equation}\label{QQ}
\begin{aligned}
\mathbf{Q}(m,n)=&\frac{2m+1}{m(m+1)}\left[\frac{m^2}{2m-1}\mathbf{P}(m,n)+\frac{1}{(2m-3)(2m+1)}\mathbf{P}(m-1,n)\right.\\
&+\left.\sum_{i=1}^{m-2}\left(\frac{1}{(2i-1)(2i+3)}\mathbf{P}(i,n)\right)-\frac{1}{3}\mathbf{P}(0,n)\right].
\end{aligned}
\end{equation}
Inserting (\ref{QQ}) into (\ref{P}) gives the following recursion,
\begin{equation}\label{nonh}
\begin{aligned}
&\mathbf{P}(m+1,n)=\left(\frac{n+1}{2n+1}\right)^2\left(\frac{2m+1}{m+1}\right)^2\mathbf{P}(m,n+1)+\left(\frac{n}{2n+1}\right)^2\left(\frac{2m+1}{m+1}\right)^2\mathbf{P}(m,n-1)\\
&\quad+\frac{2m(2m+1)(n^2+n-m^2-m)}{(m+1)^3(2m-1)(2n+1)^2}\mathbf{P}(m,n)-\left[(\frac{m}{m+1})^2-\frac{2(n^2+n-m^2-m)}{m(m+1)^3(2m-3)(2n+1)^2}\right]\mathbf{P}(m-1,n)\\
&\quad+\frac{2(n^2+n-m^2-m)(2m+1)}{m(m+1)^3(2n+1)^2}\sum_{i=1}^{m-2}\left(\frac{1}{(2i-1)(2i+3)}\mathbf{P}(i,n)\right)-\frac{4(2m+1)(n^2+n-m^2-m)}{3m(m+1)^3(2n+1)^3}.
\end{aligned}
\end{equation}
Calculating
\begin{equation*}
  \frac{m(m+1)^3(2n+1)^2}{2(n^2+n-m^2-m)(2m+1)}\mathbf{P}(m+1,n)-\frac{(m-1)m^3(2n+1)^2}{2(n^2+n-m^2+m)(2m-1)}\mathbf{P}(m,n),
\end{equation*}
by substituting the expression of recursion (\ref{nonh}) into  $\mathbf{P}(m+1,n)$ and $\mathbf{P}(m,n)$, we can obtain the recursion formula which the expression and the coefficients (\ref{coe}) are asserted in Theorem 3.3.\\
\qed\\

Here we would like to mention that the ``nonhomogeneous term'' (the last term in (\ref{nonh}) only with respect to $m$ and $n$ except $\mathbf{P}(\cdot,\cdot$)) is also disappeared with the summation about $\mathbf{P}(i,n)$!

So $\mathbf{P}(m+1,n)$ is completely finite linear dependent and this is important to get the estimate of integral of (\ref{integral}).
\subsection{Proof of lemma \ref{lp5}}
Using the identity (\ref{idd}), we can \emph{fold }the \emph{unilateral} sequence as follows
\begin{equation}\label{C}
  \mathbf{C}(m)=\beta_{m-1}\mathbf{C}(m-1)-\beta_{m-2}\mathbf{C}(m-2),
\end{equation}
where
\begin{equation}\label{Cc}
  \mathbf{C}(m)=\mathbf{c}(m+1)-\mathbf{c}(m),
\end{equation}
and
\begin{equation*}
  \beta_{m-1}=\alpha_m-1=\frac{8m^4-4m^2+1}{2(m+1)^3(2m-1)},\quad \beta_{m-2}=\alpha_{m-2}=\frac{(m-1)^3(2m+1)}{(m+1)^3(2m-1)}.
\end{equation*}
and it is easy to calculate that
\begin{equation*}
  \mathbf{C}(1)=-\frac{5}{32},\quad \mathbf{C}(2)=-\frac{5}{64}.
\end{equation*}
\begin{rem}
It is important to point out that $\mathbf{c}(m)$ is ``homogeneous'' (has no terms only with respect to $m$) and it belongs to the closure of the set $\mathscr{C}$,
\begin{equation}\label{c_m}
  \mathbf{c}(m)\in \overline{\mathscr{C}},\quad \mathscr{C}=\left\{\frac{\mathscr{P}(m)}{\mathscr{Q}(m)}\right\}
\end{equation}
where $\mathscr{P}(m)$ and $\mathscr{Q}(m)$ are polynomials with integral coefficients (with respect to $m$ ). Then it excludes the case which has the ``nonhomogeneous term'', for example
\begin{equation*}
\mathbf{\tilde{c}}(m)=\frac{m-1}{m}\mathbf{\tilde{c}}(m-1)+\frac{1}{m^2},\quad \mathbf{\tilde{c}}(1)=1
\end{equation*}
and we know $\mathbf{\tilde{c}}(m)=\frac{1+\frac{1}{2}+\ldots+\frac{1}{m}}{m}$, due to the nonhomogeneous term $\frac{1}{m^2}$.
\end{rem}
If we denote another sequence which can be \emph{folded},
\begin{equation*}
  \mathbf{D}(m)=-\frac{15}{16(m+1)(m+2)}=\frac{15}{16(m+2)}-\frac{15}{16(m+1)},
\end{equation*}
and the error sequence follows
\begin{equation*}
  \mathbf{E}(m)=\mathbf{D}(m)-\mathbf{C}(m),
\end{equation*}
then it is easy to check that, for $m\geq 3$,
\begin{equation*}
  \mathbf{D}(m)\geq \beta_{m-1}\mathbf{D}(m-1)-\beta_{m-2}\mathbf{D}(m-2),
\end{equation*}
\begin{equation*}
  \mathbf{E}(m)\geq\beta_{m-1}\mathbf{E}(m-1)-\beta_{m-2}\mathbf{E}(m-2).
\end{equation*}
It is obvious to see that
\begin{equation*}
\mathbf{E}(1)=\mathbf{E}(2)=0.
\end{equation*}
By induction, we obtain that
\begin{equation}
\mathbf{E}(m)\geq\left(\frac{m+\sqrt{2}-1}{m+\sqrt{2}}\right)^2\mathbf{E}(m-1).
\end{equation}
Then for any~ $m$,
\begin{equation*}
  \mathbf{E}(m)\geq 0,
\end{equation*}
and it follows that
\begin{equation}\label{CD}
\mathbf{C}(m)\leq \mathbf{D}(m)=-\frac{15}{16(m+1)(m+2)}<0,
\end{equation}
by the notation (\ref{Cc}), which means
\begin{equation}
0\leq \mathbf{c}(m+1)<\mathbf{c}(m)\leq \frac{1}{2}.
\end{equation}
and $\mathbf{c}(m)$ has a limitation which we will know it is $0$ in the end.

Besides, from (\ref{CD}), we can deduce that
\begin{equation*}
  \sum_{i=1}^{m}\mathbf{C}(i)\leq \sum_{i=1}^{m}\mathbf{D}(i).
\end{equation*}
\emph{unfolding }the sequence, we obtain
\begin{equation*}
 \sum_{i=1}^{m}(\mathbf{c}(i+1)-\mathbf{c}(i))\leq \sum_{i=1}^{m}\frac{15}{16}\left(\frac{1}{m+2}-\frac{1}{m+1}\right),
\end{equation*}
\begin{equation*}
\mathbf{c}(m+1)\leq\mathbf{c}(1)+\frac{15}{16}\left(\frac{1}{m+2}-\frac{1}{2}\right).
\end{equation*}
Then we have the upper bound for $\mathbf{c}(m)$,
\begin{equation}
\mathbf{c}(m+1)\leq\frac{15}{16(m+2)}+\frac{1}{32}.
\end{equation}
\qed

\section{Acknowledgements}

The authors are very grateful to the referees for their invaluable suggestions and comments.\\

\bibliographystyle{amsplain}

\end{document}